\documentclass{article}
\usepackage[T2A]{fontenc}
\usepackage[cp1251]{inputenc}
\usepackage[english,russian]{babel}
\usepackage[tbtags]{amsmath}
\usepackage{amsfonts,amssymb,mathrsfs,amscd}

\let\No=\textnumero

\usepackage[hyper]{izv2e}
\JournalName{СЕРИЯ МАТЕМАТИЧЕСКАЯ}
\numberwithin{equation}{section}
\theoremstyle{plain}
\newtheorem{theorem}{Теорема}
\newtheorem{lemma}{Лемма}[section]

\newtheorem{corol}{Следствие}
\theoremstyle{definition}
\newtheorem{problem}{Задача}

\def\klabel#1{}

\begin{document}
	 
\udk{517.984.50}

\date{}

\author{Р.\,С.~Сакс}
\address{Институт математики c ВЦ УФИЦ РАН, \\ ул.Чернышевского, 112, 450077, г. Уфа, Россия
}
\email{romen-saks@yandex.ru}

%\begin{document}

\title{Операторы ротор и градиент дивергенции в  пространствах      
		$ \mathbf{W}^{m}$ и  $\mathbf{A}^{2k}$  вихревых и потенциальных полей и в классах  $\mathbf{C}(2k, m)$}

\markboth{Р.\,С.~Сакс}{Операторы $\mathrm{ROT}$ и $\nabla\,\mathrm{DIV}$  	в пространствах Соболева}

\maketitle

\begin{fulltext}

\begin{abstract} С.Л. Соболев, 	изучая  краевые задачи для  полигармонического уравнения $\Delta^m\,u=\rho$ в пространствах ${W} _2^{(m)}(\Omega)$ c обобщённой правой частью, заложил фундамент теории этих  пространств. % \cite{sob}  $\S 9$ гл. 12.
	
	Операторы   градиент дивергенции  и  ротор ротора 
	($\nabla \text{div}$   и	 	$ \text{rot}^2$)  и их степени являются аналогами скалярного оператора 	$\Delta^m$  в   ортогональных подпространствах  $\mathcal{A}$  и $\mathcal{B}$  в  $\mathbf{L}_{2}(G)$.
	
	Они порождают  пространства $ \mathbf{A}^{2k}$ и	$\mathbf{W}^m$ потенциальных и  вихревых    полей,  а 
	собственные поля $\nabla \text{div}$   и	 	$ \text{rot}$
	задают базисы в этих пространствах. 
	
	Прямые суммы  $ \mathbf{A}^{2k}(G) \oplus \mathbf{W}^m(G)$ - новое  семейство пространств.  
	Его элементы,	 
	классы $ \mathbf{C}(2k, m)\equiv \mathbf{A}^{2k}\oplus \mathbf{W}^m$, играют роль 	${W} _2^{(m)}(G)$ 	в ограниченной области   $ G$  в 
	$\mathbb{R}^3$ с гладкой границей $\Gamma$. 
	
	Кроме того, шкалы  	$ \{\mathbf{A}^{2k}\}$  и $\{\mathbf{W}^m\}$  пространств 	$ \mathbf{A}^{2k}(G)$ и $\mathbf{W}^m(G)$ представляют   	шкалу  $\{\mathbf{H}^{n}\}$  пространств Соболева  $\mathbf{H}^{n}(G)$  в ортогональных подпространствах   $\mathbf{A}^0$ и  $ \mathbf{V}^0 $ в $\mathbf{L}_{2}(G)$, а именно, 	
	
  если поле  $\mathbf{f}\in \mathbf{H}^n(G)$ и  $n\geq 2k, n\geq m$, то его проекции   $\mathbf{f}_\mathcal{A}$ и  $\mathbf{f}_\mathbf{V}$ на подпространства 
  $\mathbf{A}^0$ и  $ \mathbf{V}^0 $   принадлежат  $\mathbf{A}^{2k}(G)$  и $\mathbf{W}^m(G) $, соответственно.  
  
  Обратно, 
   $\mathbf{A}^{2k}(G) \subset \mathbf{H}^{2k}(G)$, \quad $\mathbf{W}^m(G) \subset 
  \mathbf{H}^{m}(G)$ при   $k\geq 1, m\geq1$,  значит, 
  шкалы $ \{\mathbf{A}^{2k}\}$  и $\{\mathbf{W}^m\}$ 
     играют роль    	шкалы   $\{\mathbf{H}^{n}\}$  в    $\mathbf{A}^0$ и  $ \mathbf{V}^0 $. 
  
  В пространствах 	$ \mathbf{A}^{2k}(G)$ и $\mathbf{W}^m(G)$ 
  действуют операторы	 $\mathcal{N}_d$ и $S$, -  самосопряженные расширения 
  $\nabla \text{div}$ и  $ \text{rot}$    в  $ \mathbf{A}^0$ и  $\mathbf{V}^{0}$.  Они осущесвляют гомеоморизмы соседних пространств, а их сгепени   $\mathcal{N}^{2k}_d$ и $S^{2m}$- гомеоморизмы  пространств $ \mathbf{A}^{2k}(G)$ и $ \mathbf{A}^{-2k}(G)$ и  $\mathbf{W}^m(G)$  и  $\mathbf{W}^{-m}(G)$, соответственно.
  
Рассмотрены краевые задачи  для   $\nabla \text{div}+\lambda I$ 	 	   и	 	$ \text{rot}+\lambda I$   и других модельных операторов в  классах $  \mathbf{C}(2k, m)$.

	Приведены аргументы астрофизиков и специалистов по физике плазмы  в пользу существования безсиловых  полей (полей Бельтрами) в Природе.

\end{abstract}

 {\bf   Введение: как устроена шкала  пространств Соболева $\mathbf{H}^{m}(G)$?  }\footnote{Открылась Бездна звезд полна,
 звездам числа нет,  Бездне - дна. М.В.Ломоносов}
%   $\mathbf{A}^{2k}(G)$ и 	$\mathbf{W}^m(G)$.
 %В работе изучается шкала  
 
 Это цепь пространств  $\mathbf{H}^{m}(G)$  векторных полей в  $\mathbb{R}^3$, вложенных  подпространств и надпространств  $\mathbf{L}_{2}(G)$: 
\[\subset \mathbf{H}^{m}(G)\subset...\subset \mathbf{H}^1(G)\subset \mathbf{L}_2(G)\subset  \mathbf{H}^{-1}(G)\subset...\subset \mathbf{H}^{-m}(G) \subset  ...      \eqno{(1)}\]
 Она широко используется в теоретических и прикладных исседованиях.

С.Л. Соболев рассмотрел  её  в   книге 
{\it Введение в теорию кубатурных формул} \cite{sob}. В  $\S 9$  главы  12	
он исследовал  задачу $\pi$ и краевые задачи $D$ и $N$ для скалярного полигармонического уравнения $\Delta^m\,u=\rho$ в пространствах ${W} _2^{(m)}(\Omega)$ c правой частью -- обобщённой функцией;

 в периодическом случае (задача  $\pi$)  он доказал:  

Т е о р е м а  XII.13 (С.Л.Соболева).  {\it Оператор $\Delta^m$ переводит произвольную функцию $u$  из   $\bar {W} _2^{(m)}$ в 
	$\Delta^m\,u=\rho$ --  элемент  $\bar {L}_2^{(m)^*}$.
	Обратно, для произвольной обобщённой функции $\rho$  из
	$\bar {L}_2^{(m)^*}$ существует   функция $u\in \bar {W} _2^{(m)}$  такая, что $\Delta^m\,u=\rho$. 
	Эта функция определяется  с точностью  до  произвольного постоянного слагаемого. }

   В данной статье  представлены   две шкалы   	$ \{\mathbf{A}^{2k}\}$  и $\{\mathbf{W}^m\}$: 
\[\subset \mathbf{A}^{2k}\subset...\subset\mathbf{A}^{2}\subset  \mathbf{A}^0\subset \mathbf{A}^{-2}\subset...\subset \mathbf{A}^{-2k}\subset ...,     \eqno{(2)} \]
\[\subset \mathbf{W}^{m}\subset...\subset \mathbf{W}^1\subset \mathbf{V}^0\subset  \mathbf{W}^{-1}\subset...\subset \mathbf{W}^{-m}\subset ...,     \eqno{(3)} \]	  
введённые  автором  в работах  \cite{ds18} и  \cite{saVS20}, на базе подпространств $\mathbf{A}^0$ и $\mathbf{V}^0$  в  $\mathbf{L}_{2}(G)$.

Потенциальные (irrotational) поля $\mathbf{f}$ и соленоидальные поля $\mathbf{g}$ в  $\mathbf{L}_{2}(G)$ впервые выделил Герман Вейль в работе \cite{hw} 1941 г.
 условиями ортогональности:
\[	 (\mathbf{f},  \mathrm{rot}\,\mathbf{v})=0    \quad   \forall\, \mathbf{v}\in \mathbf{C}_0^1(G) , \quad (\mathbf{g},\nabla \psi )=0    \quad   \forall\, \psi\in C_0^1(G)  \eqno{(4)}\]

 С.Л. Соболев в статье \cite{sob54}
 %\footnote{ Об одной новой задаче	математической %физики,  Известия АН СССР (серия математическая), 
 	 1954 г.  предложил другой способ  разложения  $\mathbf{L}_{2}(\Omega)$. Область $\Omega$ у него гомеоморфна шару.

Мы используем разложение  Z.Yoshida и Y.Giga,   упростив обозначения в \cite{yogi}. Авторы  называют его разложением Вейля \cite{hw}.
 
Итак, пространство  
 ${\mathcal{{A}}}(G) =\{\nabla h, h\in H^1(G)\}$,  
а  $\mathcal{B}$ -  ортогонально  $\mathcal{A}$, 
то-есть
\[\mathbf{L}_{2}(G)=
{\mathcal{{A}}}(G)\oplus{\mathcal{{B}}}(G).  \eqno(5)\]     В свою очередь, 
\[ \mathcal{A}(G)=\mathcal{A}_{H}  \oplus \mathbf {A} ^{0}, \quad 
\mathcal{B}(G)=\mathcal{B}_{H} \oplus \mathbf {V} ^{0}.
      \eqno{(6)}\]
Здесь $\mathcal{B}_H$  и   $ \mathcal{A}_H$ - конечномерные
  подпространства размерности  $\rho$  и  $\rho_1\geq  \rho$,
  где число  $\rho$ есть  род границы $\Gamma$ области $G$ 
  (см. W. Borchers и H. Sohr  \cite{boso}).

  Пространство   $\mathcal{B}_H$ состоит из  безвихревых соленоидальных полей:
\[\mathcal{B}_H=\{\mathbf{u}\in\mathbf{L}_{2}(G):\,\mathrm{div} \mathbf{u}=0, \,\, \mathrm{rot} \mathbf{u}=0  \quad \text{в}\quad G, \quad \gamma(\mathbf{n}\cdot \mathbf{u})=0 \},  \eqno(7) \]
решений эллиптической краевой задачи. Его     базисные поля
$\mathbf{h}_j\in \mathbf{C}^{\infty}\,(\bar{G})$,    $j=1,..,\rho$;
их гладкость внутри области  $G$ впервые заметил Г.Вейль \cite{hw}. 
 Базисные поля $\mathbf{g}_l \in  \mathcal{A}_H $  также принадлежат  $\mathbf{C}^{\infty}\,(\bar{G})$,  $l=1,..,\rho_1$ (см. п.1.7). 

 Операторы	$\nabla \text{div}$ и  $ \text{rot}$  имеют самосопряженные расширения $\mathcal{N}_d$ и $S$ в   пространства  $ \mathbf{A}^0$ и  $\mathbf{V}^{0}$ потенциальных  и  вихревых полей.
  Оператор  $S$  опредедён  Z.Yoshida and Y.Giga в \cite{ yogi},  а оператор $\mathcal{N}_d$ -  автором в  \cite{ds18}.

Собственные поля операторов  $\mathcal{N}_d$ и $S$ 
в    $ \mathbf{A}^0$ и  $\mathbf{V}^{0}$
задают   ортонормированные базисы    $\{ \mathbf{q}_{j}(\mathbf{x})\}$ и  
$\{\mathbf{q}^{+}_{j}(\mathbf{x}), \mathbf{q}_{j}^{-}(\mathbf{x})\}$:
\[-\nabla\text{div}\mathbf{q}_{j}(\mathbf{x})=\nu^2_j \mathbf{q}_{j}(\mathbf{x}), \quad
\mathbf{x}\in G, \quad \gamma \mathbf{n}\cdot \mathbf{q}_{j}=0,
\quad \|\mathbf{q}_{j}\|=1,
\eqno{(8)}\]
 %в  $\mathbf{V}^{0}$базис ,$: 
\[\mathrm{rot}
\mathbf{q}_{j}^{\pm}=\pm\lambda_j\, \mathbf{q}_{j}^{\pm}, \quad   \gamma\mathbf{n}\cdot\mathbf{q}_{j}^{\pm}=0,   \quad  \|\mathbf{q}^{\pm}_{j}\|=1  \quad   j=1, 2, ....     \eqno{(9)} \]
Доказано, что $\mathbf{q}_{j}$ и  $\mathbf{q}^{\pm}_{j}\in \mathbf{C}^\infty(\bar{G})$ (см. п. 1.6).

Элементы    $\mathbf{f}_\mathcal{A}\in \mathbf{A}^0$   и  $\mathbf{f}_\mathbf{V}\in \mathbf{V}^0 $   
представляются   рядами Фурье: 
\[\mathbf{f}_\mathcal{A}=
\sum_{j=1}^{\infty} (\mathbf{f},{\mathbf{q}}_{j}) \mathbf{q}_{j }(\mathbf{x}), \quad 	\mathbf{f}_\mathbf{V}=\sum_{j=1}^{\infty}\left[(\mathbf{f},{\mathbf{q}}_{j}^{+})
\mathbf{q}_{j}^{+}(\mathbf{x})+
(\mathbf{f},{\mathbf{q}}_{j}^{-})
\mathbf{q}_{j}^{-}(\mathbf{x})\right].     \eqno{(10)}\]
а  операторы    $\mathcal{N}_d$,   $S$  и их обратные --
преобразованиями рядов:
\[  \mathcal{N}_d \mathbf{f}_\mathcal{A}=-\sum_{j=1}^{\infty} \nu^2_j  (\mathbf{f},{\mathbf{q}}_{j}) \mathbf{q}_{j }, \quad   
	S\mathbf{f}_\mathbf{V}=\sum_{j=1}^{\infty}\lambda_{j}\left[(\mathbf{f},{\mathbf{q}}_{j}^{+})
\mathbf{q}_{j}^{+}-(\mathbf{f},{\mathbf{q}}_{j}^{-})
\mathbf{q}_{j}^{-}\right],  \eqno{(11)}  \]
\[  \mathcal{N}_d  ^{-1}\mathbf{f}_\mathcal{A}=-
\sum_{j=1}^{\infty} \nu^{-2}_j  (\mathbf{f},{\mathbf{q}}_{j}) \mathbf{q}_{j }, \quad    	S^{-1}\mathbf{f}_\mathbf{V}=\sum_{j=1}^{\infty}\lambda_{j}^{-1}\left[(\mathbf{f},{\mathbf{q}}_{j}^{+}) \mathbf{q}_{j}^{+}-
(\mathbf{f},{\mathbf{q}}_{j}^{-}) \mathbf{q}_{j}^{-} \right].     \]
%Отметим, что принадлежности  $\mathbf{f}_\mathcal{A}\in %\mathbf{A}^{2}(G)$ и  $\mathbf{f}_\mathbf{V}\in %\mathbf{W}^m(G)$ %$$  и $ $,

Принадлежность  элементов    $\mathbf{f}_\mathcal{A}$ к $\mathbf{A}^0$   и  $\mathbf{f}_\mathbf{V}$   к 
$ \mathbf{V}^0 $   
определяется сходимостью  рядов: 
\[ \|\mathbf{f}_\mathcal{A}\|^2=
\sum_{j=1}^{\infty} (\mathbf{f},{\mathbf{q}}_{j})^2 , \quad  \|\mathbf{f}_\mathbf{V}\|^2=\sum_{j=1}^{\infty}\left[(\mathbf{f},{\mathbf{q}}_{j}^{+})^2
+(\mathbf{f},{\mathbf{q}}_{j}^{-})^2\right].     \eqno{(12)}\]
Затем,   $\mathbf{f}_\mathcal{A}\in \mathbf{A}^{2}$   и  $\mathbf{f}_\mathbf{V}\in \mathbf{W}^1$ тогда и только тогда, когда    
\[ \|\mathbf{f}_\mathcal{A}\|<\infty,  \quad \|  \mathcal{N}_d \mathbf{f}_\mathcal{A}\|<\infty  \quad \text{ и } \quad   \|\mathbf{f}_\mathbf{V}\|<\infty,  \quad \|S\mathbf{f}_\mathbf{V}\|<\infty ;    \eqno{(13)}\]
и так далее  (см. \cite{ds18,saVS20}).

Имеются   {\bf соответствия } между шкалами  (1) и  (2) , (3).   

Если поле  $\mathbf{f}\in \mathbf{H}^n(G)$ и  $n\geq 2k, n\geq m$, то его проекции   $\mathbf{f}_\mathcal{A}$ и  $\mathbf{f}_\mathbf{V}$ на подпространства 
$\mathbf{A}^0$ и  $ \mathbf{V}^0 $   принадлежат  $\mathbf{A}^{2k}(G)$  и $\mathbf{W}^m(G) $, соответственно.

  Обратно, пространства
$ \mathbf{W}^1\subset \mathbf{H}^{1}$,  $\mathbf{A}^{2}\subset \mathbf{H}^{2}$ согласно  оценкам  (1.12) и (1.13) при $s=0$.
По индукции, $\mathbf{A}^{2k}(G) \subset \mathbf{H}^{2k}(G)$, \quad $\mathbf{W}^m(G) \subset 
 \mathbf{H}^{m}(G)$ при   $k\geq 1, m\geq1$, соответственно.  
 
 Значит,    {\bf  шкалы (2) и (3)  
   представляют шкалу (1) 	  }
  в ортогональных подпространствах	$\mathbf{A}^0$ и $\mathbf{V}^0$. %\cite{saVS20}.

{\bf  Свойства } шкал (2) и (3).
    
1) Операторы $\mathcal{N}_d$,  $\mathcal{N}^{-1}_d$ и $S$,  $S^{-1}$  отображают  соседние пространства Гильберта в шкалах (2) и (3)    взаимно однозначно и непрерывно (Лемма 2 и Теорема 3 в \cite{saVS20}).     
Изобразим это так: 
\[\mathcal{N}_d:\mathbf{A}^{2k}(G)\to \mathbf{A}^{2(k-1)}(G),  \quad \mathcal{N}^{-1}_d:\mathbf{A}^{2
	(k-1)}(G)\to \mathbf{A}^{(2k)}(G),  \quad k>0,     \eqno{(14) } \]
\[ S: \mathbf{W}^m(G)\to \mathbf{W}^{{m-1}}(G),     \quad                          S^{-1}: \mathbf{W}^{m-1}(G)\to \mathbf{W}^{m}(G),   \quad m>0.     \eqno{(15)} \] 
 Значит, 
\[\mathcal{N}_d: \mathbf{A}^{2k}\to...\to\mathbf{A}^{2}\to  \mathbf{A}^0,   \quad \mathcal{N}^{-1}_d:\mathbf{A}^0\to \mathbf{A}^{2}\to...\to\mathbf{A}^{2k}\to... ,     \eqno{(16)} \]
\[ S: \mathbf{W}^{m}\to...\to \mathbf{W}^1\to \mathbf{V}^0,  \quad  S^{-1}:\mathbf{V}^0\to  \mathbf{W}^{1}\to ...\to  \mathbf{W}^{m}\to ....   \eqno{(17)} \]

Таким образом,  пространства $\mathbf{A}^{2k}(G)$ и 	$\mathbf{W}^m(G)$  разных порядков  $k$ и   $m$  топологически устроены одинаково.

2) Спектры операторов 	$-\nabla \text{div}$  и $ \text{rot}$
точечные и действительные \cite{yogi,saVS20} .

3) Если      область       $G=B$ есть шар  радиуса $R$, то        собственные значения  ротора   равны  $\lambda_{n,m}=\pm \rho_{n,m}/R$,  где числа $\pm \rho_{n,m}$ - нули элементарных  функций  
\begin{equation*}\label{bes  2}
	\psi_n(z)
	= (-z)^n\left(\frac d{zdz}\right)^n\left(\frac{\sin z}z\right),  \quad  m,n\in {\mathbb {N}}. \eqno{(18)}
\end{equation*}
Функции  $\psi_n(z)\equiv J_{n+1/2}(z)$, цилиндрические  функции   полуцелого порядка, - элементарны,  это заметил ещё Леонард Эйлер (см. \cite{vla} $\S23$).

Кратность  собственного значения  $\lambda^{\pm}_{n,m}$ равна $2n+1$. 

Собственные значения  оператора $-\nabla\mathrm{div}$ равны $\nu_{n,m}^2$, где $\nu_{n,m}=\alpha_{n,m}/R$,  а  числа  $\alpha_{n,m}$ - нули производных  $\psi'_n(r)$,   $n \geq 0,\, m\in {\mathbb {N}}$ ;     	
кратность  собственных значений  $\nu^2_{n,m}$    равна $2n+1$.	 	

Приводятся формулы  (1.28) и (1.35) собственных полей
$ \mathbf{u}_{\kappa}^{\pm }$  и  $ \mathbf{v}_{\kappa} $  ротора  и $-\nabla\mathrm{div}$ в сферических координатах,  вычисленные  автором в
 \cite{saUMJ13,saVSTU}.

В этом случае, пространства $\mathbf{A}^{2k}(B)$ и $\mathbf{W}^{{m}}(B)$ определяются явно и можно использовать вычислительные средства.

4) Имеется симметрия между пространствами положительных и отрицательных порядков:
\[\mathcal{N}^{2k}_d:\mathbf{A}^{2k}(G)\to \mathbf{A}^{-2k}(G),  \quad \mathcal{N}^{-2k}_d:\mathbf{A}^{-2k}(G)\to \mathbf{A}^{2k}(G),  \quad k>0,     \eqno{(19) } \]
\[ S^{2m}: \mathbf{W}^{m}(G)\to \mathbf{W}^{{-m}}(G),     \quad                          S^{-2m}: \mathbf{W}^{-m}(G)\to \mathbf{W}^{{m}}(G),   \quad m>0.     \eqno{(20)} \] 
Точнее, если  $\mathbf{u}$ произвольный элемент из $\mathbf{A}^{2k}$,  то $\mathbf{v}=(\nabla \mathrm{div})^{2k}\,\mathbf{u}$-функционал из $(\mathbf{A}^{2k}_0)^*$, 
определяемый интегралом:
\[\int_G
\mathbf{v}\cdot (\mathbf{w}_{\eta})\, d\mathbf{x}=	
((\nabla \mathrm{div})^{2k}\,\mathbf{u},\,\mathbf{w}_{\eta})\equiv 
((\nabla \mathrm{div})^{k}\,\mathbf{u},(\nabla \mathrm{div})^{k}\,\mathbf{w}_{\eta})`.   \eqno{(21)}\]
 где 
$\mathbf{w}_{\eta}$ -- средняя вектор-функция для поля $\mathbf{w}\in \mathbf{A}^{2k}$,  поле
$\mathbf{w}_{\eta}\in\mathbf{A}^{2k}_0$ и равно нулю вблизи границы области $G$.

Пространство $(\mathbf{A}^{2k}_0)^*$ отождествляется с
пространством $\mathbf{A}^{-2k}$.

Если  пространство $ \mathcal{A}_H$ не пусто, то функционал $\mathbf{v}$  на $ \mathcal{A}_H$ равен нулю.

%Доказано, что

Соотношения (10), (11)  обосновываются утверждениями:

Т е о р е м а  4.
	{\it При заданном $\mathbf{v}$   в объединении  $ \mathcal{A}^{-2k}$ и $k\geq 1$
	уравнение $(\nabla \mathrm{div})^{2k}\,\mathbf{u}=\mathbf{v}$  разрешимо в
	пространстве $\mathbf{A}^{2k}$ тогда и только тогда, когда
	$\mathbf{v}\in\mathbf{A}^{-2k}$.\quad
	Его решение   $\mathbf{u}= \mathcal{N}_d^{-2k}\mathbf{v}$
	в  фактор-пространстве	$\mathcal{A}/\mathcal{A}_H$   определяется однозначно.}

Т е о р е м а  5.
{\it 
	При заданном $\mathbf{v}$  в объединении  $\mathbf{W}^{-m}$
	уравнение $\mathrm{rot}^{2m}\,\mathbf{u}=\mathbf{v}$  разрешимо в
	пространстве $\mathbf{W}^m(G)$ тогда и только тогда, когда
	$\mathbf{v}\in \mathbf{W}^{-m}(G)$.\,
	Его решение $\mathbf{u}=S^{-2m}\mathbf{v}$ в   фактор-пространстве	
	$ \mathcal{B}(G)/\mathcal{B}_H(G)$  определяется 
	однозначно}.

Теоремы 4 и 5 доказаны в $\S 2, 3$. Их  доказательства  прояснились  при чтении $\S 9$ главы 12 основополагающей  книги  С.Л.Соболева \cite{sob}.

В Теореме 6 п.3.5, доказанной  в \cite{saVS20},  утверждается, что операторы 
 $(S+\lambda\,I):\mathbf{W}^{m+1} \rightarrow \mathbf{W}^{m}$ при $m\geq 0$  и   $(S+\lambda\,I)^{-1}:\mathbf{W}^{m} \to \mathbf{W}^{m+1}$  разрешимы по Фредгольму, а при 	при $\lambda \overline{\in} Sp (\mathrm{rot})$ они однозначно обратимы, приводятся оценки норм.
 
В Теореме 7 п.3.6  приводятся соотношения между пространствами  $ \mathbf{W}^{k}(\Omega)$,    $ \mathbf{H}^{k}(\Omega)$	и    $ \mathbf{C}^{k-2}(\bar{ \Omega})$ в области $\Omega$,  гомеоморфной шару,  доказаные в \cite{saVS20}.  

В $\S 4$ мы рассмотрели   новые   классы пространств     $ \mathbf{C}(2k,m)\equiv\mathbf{A}^{2k} \oplus \mathbf{W}^m$   в  области $\Omega$,  гомеоморфной шару; где числа $k, m $  -целые.
 Они образуют сеть или дву-параметрическое  семейство пространств 
и принадлежат     $\mathbf{L}_{2}(G)$, если  $k\geq 0$ и  $m\geq 0$.
Класс  $  \mathbf{C}(2k, 2k)$ совпадает с пространством Соболева    $\mathbf{H}^{2k}(G)$.

Ставятся краевые задачи  для $(\nabla \mathrm {div})^p +\lambda I $ и $ (\text{rot})^p+\lambda I$  при  $p=1,2, ...$  и других модельных операторов в  пространствах  $ \mathbf{C}(2k, m)$.  

Открылось широкое поле задач.  Мы рассмотрели простейшие из них, указав путь их решения. 
 В статье  \cite{saVS20} доказано

{\it	При $\lambda\neq Sp (rot)$ единственное решение задачи 1.1  имеет вид:
		\[\mathbf{u}={\mathbf{u}_{\mathcal{A}}}+{\mathbf{u}_{\mathbf{V}}}, \quad \text{где}\quad 
		{\mathbf{u}_{\mathcal{A}}}={\lambda}^{-1}\mathbf{f}_\mathcal{A}, \quad 
		\mathbf{u}_\mathbf{V}=(S+\lambda \, I)^{-1} \mathbf{f}_\mathbf{V}.		\eqno{(22)}\]
	
Это решение $\mathbf{u}\in  \mathbf{C}(2k, m+1) $ при	$\mathbf{f}\in  \mathbf{C}(2k, m)$, где  числа  $k, m > 0 $  -целые.
Если	же $\mathbf{f}\in \mathbf{C}^{\infty}_0(\Omega)$, то поле
$\mathbf{u}\in  \mathbf{C}^{\infty}(\overline{\Omega})$ есть классическое решение задачи.}(Теорема  8).

Следствие. 
		{\it	Если область $\Omega=B$,  %есть шар, 
	$\psi_n(\lambda\,R)\neq 0$   $\forall\, n\in {\mathbb {N}}$,   числа $k, m>0$ целые, а поле  $\mathbf{f}\in  \mathbf{A}^{2k}(B)\oplus  \mathbf{W}^m(B)$,   то решение  задачи 1.1  существует, единственно и принадлежит классу    $ \mathbf{A}^{2k}(B)\oplus  \mathbf{W}^{m+1}(B)$.  }

Лемма 4.1.
{\it  		При $\lambda \overline{\in} Sp (\mathrm{rot})$ и $m>0$ оператор  $\mathrm{rot}+\lambda I$  (и обратный)
	отображает класс  $ \mathbf{C}(2k, m+1)$	на  класс  $ \mathbf{C}(2k, m)$	взаимно	однозначно и непрерывно. }

Аналогичные результаты (Терема 9, Следствие и Лемма 4.1) имеют место для задачи 2.1,  а Терема 10 - для задачи 3.

В $\S 5$ мы приводим  аргументы астрофизиков и иследователей физики плазмы в пользу  того,  что поля Бельтрами {\bf существуют в Природе} в микро- и макро-мире! 
В докладе НАН  \cite{chawo} \  "О бессиловых магнитных полях" \newline S. Chandrasekhar and L.Woltjer отмечают, что
 поле в Крабовидной туманности имеет необычайно ровный узор.  Согласно   \cite{cdtgt} это поле  $u_{1,1,0}^{+}(\mathbf {x})$.

	Профессор Исламов Г.Г.\cite{gais}, используя  формулы из  \cite{saVSTU} и  программу  Wolfram Mathematica,  осуществил визуализацию линий тока поля   $u_{1,1,0}^{+}(\mathbf {x})$ ротора  в шаре радиуса 1 со значением $\rho =\rho_{1,1}$   (см. п.1.11).
 
 Его изящный рисунок представлен в конце статьи \cite{saVS20}

	J. Cantarella,  D. De Turck, H. Gluck and  M.Teitel 
	 ("Физика плазмы")   исследовали собственные поля ротора  в шаре  радиуса $b$  и в шаровом слое.

Уравнение (1.31) на функцию  $v=r\,u_r$  при минимальном   собственном значении $\lambda_{1,1}=\rho_{1,1}/b>0$ автор обнаружил  в их статье \cite{cdtgt}.
Они приводят  также   соответствующую	 $\lambda_{1,1}$  формулу  собственного поля ротора в шаре  (см. Theorem A). К сожалению,  с опечаткой... ( Может это шифр-код?). 

Исправив её,   мы пришли к    компонентам   (1.37)  поля  $\mathbf{u}^{+}_{(1,1,0)}(\mathbf{x})$. 

Вот, что   пишут  H.Qin, W.Liu, H.Li, and J.Squire в отчёте  PPPL-4823  \footnote{ Prinston Plasma Physics Laboratory  "Woltjer-Taylor state without Taylor's conjecture-plasma relaxation  at all wavelengths"\, Okt. 2012. } {\it
	о релаксации плазмы:}
"в астрофизике  и в  физике  рлазмы было открыто, что плазма релаксирует к  состоянию Woltjier- Taylor'а  \cite{wol},  \cite{tay},  которое определяется уравнениями  Бельтрами $\nabla\, \times\,\textbf{B}= \alpha\, \textbf{B}$ с постоянной $\alpha$"\

*

	%Библиография: 36 названий.
% \end{document}

\section { Основные подпространства $\mathbf{L}_{2}(G)$}

\subsection {Шкала  пространств Соболева} Мы рассматриваем линейные
пространства над полем $\mathbb{R}$ действительных чисел. Через
$\mathbf{L}_{2}(G)$ обозначаем пространство Лебега вектор-функций  (полей),
квадратично интегрируемых в $G$ с внутренним произведением
 $(\mathbf {u},\mathbf {v})= \int_G \mathbf
{u}\cdot \mathbf {v}\,d \mathbf {x}$ и нормой
$\|\mathbf{u}\|= (\mathbf {u},\mathbf {u})^{1/2}$.

 Пространство Соболева, состоящее из полей,
 принадлежащих $\mathbf{L}_{2}(G)$ вместе с обобщенными производными
 до порядка $ m> 0$, обозначается через
$\mathbf{H}^{m}(G)$, $\|\mathbf {f}\|_m$ -норма его элемента
$\mathbf {f}$; $\mathbf{H}^{0}(G)\equiv\mathbf{L}_{2}(G)$. 

 $\mathbf{H}^m(G)$ - гильбертово пространство со 
  скалярным произведением:
\[	 (\mathbf{f},\mathbf{g})_m=(\mathbf{f},\mathbf{g})+
\int_{G} \sum_{|\alpha|=m}\frac{m!}{\alpha
	!}\partial^{\alpha}\mathbf{f}\cdot\partial^{\alpha}\mathbf{g}
d \mathbf{x},  \quad   \|\mathbf {f}\|_m^2= (\mathbf{f},\mathbf{f})_m.     \eqno{(1.1)}\] 	

Замыкание в норме $\mathbf{H}^{m}(G)$ множества $[\mathcal{C}^{\infty}_0(G)]^3$
обозначается через $\mathbf{H}^{m}_0(G)$.

 Пространство Соболева
отрицательного порядка $\mathbf{H}^{-m}(G)$ двойственно к
$\mathbf{H}^{m}_0(G)$

 На лекции    в  НГУ  в 1962  году  С.Л.Соболев рисовал всю цепь  вложенных пространств: 
 \[\subset \mathbf{H}^{m}\subset...\subset \mathbf{H}^1\subset \mathbf{L}_2\subset  \mathbf{H}^{-1}\subset...\subset \mathbf{H}^{-m} \subset       \eqno{(1.2)}\]
Он обозначал их  $W_2^{(m)}(G)$ в  \cite{sob}  $\S 9$ гл. 12. 

 Их обозначают также $\mathbf{H}^{m}(G)$ (см.  В.П.Михайлов \cite{mi} $\S 4$ гл. 3).
 
 В области $G$ с гладкой границей $\Gamma$ в каждой точке $y\in\Gamma$
 определена нормаль $\mathbf {n}(y)$ к $\Gamma$.
Поле $\mathbf {u}$ из $\mathbf{H}^{m+1}(G)$ имеет след $
\gamma(\mathbf {n}\cdot\mathbf {u})$ на $\Gamma$ его нормальной
компоненты, который принадлежит пространству Соболева-Слободецкого
$\mathbf{H}^{m+1/2}(G)$, $|\gamma(\mathbf {n}\cdot\mathbf
{u})|_{m+1/2}$- его норма.

\subsection {Ортогональные  подпространства
	$\mathcal{A}$ и $\mathcal{B}$ в $\mathbf{L}_{2}(G)$}
Пусть $h$- функция  из ${H}^{1}(G)$, а $\mathbf{u}=\nabla h$ - ее
градиент. 
По определению ${\mathcal{{A}}}(G) =\{\nabla h, h\in H^1(G)\}$,  
 а  $\mathcal{B}$ - ортогональное дополнение  $\mathcal{A}$ в пространстве  $\mathbf{L}_{2}(G)$.  

 Из соотношений ортогональности  $(\mathbf {u},\nabla h)=0$ для любой   $ h\in H^1(G)$ при $\mathbf{u}\in\mathbf{H}^{1}(G)$ вытекает, что $ \mathrm{div} \mathbf{u}=0$  в  $G$, \,$\gamma(\mathbf{n}\cdot \mathbf{u})=0 $.  Поэтому  ${\mathcal{{B}}}(G)$ обозначают ещё так:  $\mathcal{B}(G)=\{\mathbf{u}\in\mathbf{L}_{2} (G): \mathrm{div} \mathbf{u}=0 \quad \text{в}\quad G, \,\gamma(\mathbf{n}\cdot \mathbf{u})=0 \}$.  
 \footnote{Если   $\mathbf{u}$  и  $ \mathrm{div}\mathbf{u}\in\mathbf{L}_{2} (G)$,  то след $\,\gamma(\mathbf{n}\cdot \mathbf{u})$ существует \cite{saVS20}.} 
 Итак, 
\[\mathbf{L}_{2}(G)=
{\mathcal{{A}}}(G)\oplus{\mathcal{{B}}}(G).  \eqno{(1.3)}\]
{\it Замечание.} {\small Это разложение содержится в статье   Z.Yoshida и Y.Giga    \cite{yogi}}.   Авторы  называют его разложением Вейля \cite{hw}, а  ${\mathcal{{B}}}(G)$  обозначают как ${L}_{\sigma}^2(G)$.

 Если граница $\Gamma$ имеет положительный род $\rho$,  то 
   $ \mathcal{A}$  содержит  подпространство 
   \[\mathcal{A}_H=\{\mathbf{v
	}\in\mathbf{L}_{2}(G):\,\nabla\text{div}\mathbf{v}=0, \quad
 \text{rot}\,\mathbf{v}=0 \quad \text{в}\quad G, \quad \gamma (\mathbf{n}\cdot \mathbf{v})=0 \},   \eqno{(1.4)}\]  
а  $ \mathcal{B}$ содержит  подпространство безвихревых соленоидальных полей
\begin{equation*}\klabel{bh 1} \mathcal{B}_H=\{\mathbf{u}\in\mathbf{L}_{2}(G):\,\mathrm{div} \mathbf{u}=0, \,\, \mathrm{rot} \mathbf{u}=0  \quad \text{в}\quad G, \quad \gamma(\mathbf{n}\cdot \mathbf{u})=0 \}.  \eqno(1.5) \end{equation*}
  Размерность $\mathcal{B}_H$ равна $\rho$ \cite{boso} и
его     базисные поля
 $\mathbf{h}_j\in \mathbf{C}^{\infty}\,(\bar{G})$,  $j=1,..,\rho$ \cite{hw}. 
 Размерность $\mathcal{A}_H$  не меньше  $\rho$, так как 
 $\mathcal{B}_H\subset \mathcal{A}_H$.  Его     базисные поля
$\mathbf{g}_l \in \mathbf{C}^{\infty}\,(\bar{G})$,  $l=1,..,\rho_1\geq\rho$ (см. п. 1.7).

Отметим, что    у сферы
 размерность $\rho =\dim\mathcal{B}_H =0$,  у тора $\rho=1$  и
 $\rho\geq 1$  у сферы с ручками (числом  $\rho$)  . 

 Ортогональное дополнение   в $\mathcal{A}$ к   $\mathcal{A}_H$  
 обозначается    $\mathbf {A}^{0} (G)$.
 
 Ортогональное дополнение   в $\mathcal{B}$  к   $\mathcal{B}_H$ обозначается  $\mathbf {V}^{0} (G)$   и называется классом {\it вихревых}\,   полей   \cite{saVSTU}. Так что   
\begin{equation*}\klabel{bhr 1}
\mathcal{A}(G)=\mathcal{A}_{H} (G) \oplus \mathbf {A} ^{0} (G), \quad 
\mathcal{B}(G)=\mathcal{B}_{H} (G) \oplus \mathbf {V} ^{0} (G). 
 \eqno(1.6) \end{equation*}
{\it  В шаре   $B$,   множества $\mathcal{A}_H$ и $\mathcal{B}_H$   пусты  и  $\mathbf {A}^{0} (B)= \mathcal{A}(B)$,  а   $\mathbf {V}^{0} (B)=\mathcal{B}(B)$.}

{\it Замечание. }  {\small  С.Л.Соболев\,\cite{sob54}, \,О.А.Ладыженская\,\cite{lad}, К.Фридрихс \,\cite{fri},  Э.Б.Быховский и Н.В. Смирнов \, \cite{bs}  также	приводят  разложения $\mathbf{L}_{2}(G)$ на   	ортогональные подпространства. 
Причём,  С.Л.~Соболев  	предполагает, что область  (он  обозначает её $\Omega$)  гомеоморфна шару.    Z.Yoshida и Y.Giga отмечают в  \cite{yogi},  что разложение $ \mathcal{B}(G)$  (1.6) содержится в книге C.B. Morrey \cite{mo}.}

  Мы   будем придерживаться разложений (1.3), (1.6). 
  
    \subsection{  Операторы   $\nabla\mathrm{div}$   и   $\mathrm{rot}$ - проекторы $\mathbf{L}_{2}(G)$ на 	$\mathcal{A}$ и $\mathcal{B}$} Операторы градиент,  ротор (вихрь) и дивергенция определяются в трехмерном векторном анализе, например,  в курсе В.А.Зорича \, \cite{zo}. Им соответствует оператор $d$ внешнего
дифференцирования на формах $\omega^k$ степени $k=0,1$ и 2
Соотношения $dd\omega^k=0$ при $k=0,1$ имеют вид $\mathrm{rot}\,\nabla h=0$ и $\mathrm{div}\, \mathrm{rot} \mathbf{u}=0$ для гладких функций $h$ и  $ \mathbf{u}$.
Следовательно,   операторы   $\nabla\mathrm{div}$ и  $\mathrm{rot}$ аннулируют друг друга: 
\[\nabla\mathrm{div}\, \mathrm{rot}\, \mathbf{u}=0, \quad     \mathrm{rot}\,\nabla\mathrm{div} \mathbf{u}=0.     \eqno{(1.7)}\]    
  Оператор Лапласа
выражается через них  и скалярный оператор $\Delta_c$: 
\begin{equation*}\klabel{dd 1}\mathrm{\Delta} \mathbf {v}
\equiv \nabla \mathrm{div}\,\mathbf {v}
-(\mathrm{rot})^2\, \mathbf {v}= \Delta_c I_3\,\mathbf {v} , \quad  %\text{где}
 \,\,  \mathbf{v}=(v_1, v_2, v_3),  \quad  \Delta_c v_j\equiv  \mathrm{div} \nabla {v_j}\eqno(1.8) \end{equation*}% \quad  
 Оператор Лапласа эллиптичен  \cite{vol, so71, sa75},  а операторы $\mathrm{rot}$ и $\nabla \mathrm{div}$ не являются эллиптическими.  Они
вырождены, причем $\mathrm{rot}\, \mathbf {u}=0$ при $\mathbf{u}\in
\mathcal{A}$,  а $\nabla\mathrm div \mathbf {v}=0$ при $\mathbf{v}\in
\mathcal{B}$ в смысле  $\mathbf{L}_{2}(G)$ \cite{hw}.   
Поэтому 
\[\Delta \mathbf{v}\equiv  \nabla \mathrm{div} \mathbf{v}   \quad  \text{при}    \quad  \mathbf{v}\in \mathcal{A},  \quad  \Delta \mathbf{u} \equiv -\mathrm{rot}\, \mathrm{rot}\mathbf{u}  \quad  \text{при}  \quad      \mathbf{u}\in \mathcal{B}.   \eqno(1.9)\]

\subsection{Краевые   задачи для $\text{rot}$    и  $ \nabla\text{div}$ в  пространствах Соболева} 
  В классе Б.Вайнберга и В.Грушина  равномерно неэллиптических псевдодифференциальных операторов  \cite{vagr}, автор выделил в  \cite{saVS20}  подкласс [REES p] обобщённо эллиптических дифференциальных операторов   и доказал, что  {\it операторы $ \text{rot}+\lambda \, I$    и  $ \nabla\text{div}+\lambda \,I$   первого и второго порядков  при $\lambda\neq 0$ принадлежат  классу} [REES 1]. В пространствах Соболева $\mathbf{H}^{s}(G)$  изучены     краевые   задачи.  Им  соответствуют операторы $\mathbb{A}$ и $\mathbb{B}$, которые расширятся до эллиптических по В.Солонникову    переопределённых операторов $\mathbb{A}_R$ и $\mathbb{B}_R$,  ограниченных в  пространствах $\mathbf{H}^{s}(G)$ при целом $s\geq 0$:
\begin{equation*}\klabel{op 1r} \mathbb{A}_R\mathbf{u}\equiv\left( \begin{matrix}
\mathrm{rot} +\lambda I \\
\lambda \, \mathrm{div}\\
\gamma
\mathbf{n}\cdot \end{matrix}\right)\mathbf{u}:
\mathbf{H}^{s+1}(G)\rightarrow\left(
\begin{matrix}\mathbf{H}^{s}(G)\\
H^s(G)\\
H^{s+1/2}(\Gamma)\end{matrix}\right),    \eqno(1.10)
\end{equation*}
\begin{equation*}\klabel{op 1r} \mathbb{B}_R\mathbf{u}\equiv\left( \begin{matrix}
\nabla\,\mathrm{div} +\lambda I \\
\lambda \, \mathrm{rot}\\
\gamma
\mathbf{n}\cdot \end{matrix}\right)\mathbf{u}:
\mathbf{H}^{s+2}(G)\rightarrow\left(
\begin{matrix}\mathbf{H}^{s}(G)\\
\mathbf{H}^{s+1}(G)\\
H^{s+3/2}(\Gamma)\end{matrix}\right).               \eqno(1.11)
\end{equation*}

Из  Теоремы 1.1  В.Солонникова \cite{so71}  в работе \cite{saVS20} получена:
 \begin{theorem} \klabel{rot   1}
	Оператор $\mathbb{A}_R$	имеет левый регуляризатор.  
	Его ядро конечномерно и для любых $ \mathbf{u}\in  \mathbf{H}^{s+1}(G)$ и $ \lambda \neq 0 $ ( с  постоянной  $C_s =C_s(\lambda )>0$,  зависящей только от $s, \lambda $)	выполняется оценка: 
		\begin{equation*}
	\klabel{arot__1_} C_s\|\mathbf{u}\|_{s+1}
	\leq\|\mathrm{rot} \mathbf{u}\|_{s}+
	|\lambda|\|\mathrm{div} \mathbf{u}\|_{s}+
	|\gamma({\mathbf{n}}\cdot\mathbf{u})|_{s+1/2}+
	\|\mathbf{u}\|_{s}.\eqno{(1.12)}
	\end{equation*}
\end{theorem} 
\begin{theorem} \klabel{Nd_1}
	Оператор $\mathbb{B}_R$ 	имеет левый регуляризатор.
	Его ядро конечномерно и   для любых $\mathbf{v}\in  \mathbf{H}^{s+2}(G)$ и $\lambda \neq 0 $ ( с  постоянной  $C_s =C_s(\lambda )>0$,  зависящей только от $s, \lambda $)	выполняется оценка:
	\begin{equation*}
	\klabel{ond__2_} C_s\|\mathbf{v}\|_{s+2}
	\leq|\lambda|\|\mathrm{rot} \mathbf{v}\|_{s+1}+
	\|\nabla\mathrm{div} \mathbf{v}\|_{s}+
	|\gamma({\mathbf{n}}\cdot\mathbf{v})|_{s+3/2}+ \|\mathbf{v}\|_{s}.\eqno{(1.13)}
	\end{equation*}
\end{theorem}

 Топологических ограничений на область нет, предполагается ее
связность, ограниченность и гладкость границы.
Оценка (1.12) известна давно (см. \cite{fri, yogi}). 
 Здесь мы показываем, что  для  операторов класса  [REES p]  аналогичные оценки
можно получать из Теоремы В.Солонникова. 

 Формулы $\mathbf{u}\cdot\nabla h+ h\mathrm{div}\mathbf{u}=\mathrm{div}(h \mathbf{u}) $, \quad 
 $\mathbf{u}\cdot\mathrm{rot} \mathbf{v}- \mathrm{rot}
 \mathbf{u}\cdot\mathbf{v}=\mathrm{div}[\mathbf{v},\mathbf{u}]$, где
 $\mathbf{u}\cdot\mathbf{v}$ и 
 $[\mathbf{v},\mathbf{u}]$ -скалярное и  векторное произведения в   $R^3$, и интегрирование по  $G$ используются при определении операторов $\nabla\,\mathrm{div}$ и $\mathrm{rot}$ в $\mathbf{L}_{2}(G)$.  
  Интегрируя и применяя формулу Гаусса-Остроградского, имеем
 \[ \int_G [\mathrm{rot} \mathbf{u}\cdot \mathbf{v}
 -\mathbf{u}\cdot \mathrm{rot} \mathbf{v}]d\mathbf{x} =  \int_{\Gamma}   \mathbf{n}\cdot [\mathbf{v},\mathbf{u}]d\mathbf{S}.  \eqno{(1.14)}\] 
 \[ \int_G [\nabla \mathrm{div} \mathbf{u}\cdot \mathbf{v}
 -\mathbf{u}\cdot \nabla\mathrm{div} \mathbf{v} ]d\mathbf{x} =  \int_{\Gamma} [ (\mathbf{n}\cdot \mathbf{v})\mathrm{div} \mathbf{u}+  (\mathbf{n}\cdot \mathbf{u})\mathrm{div} \mathbf{v}]d\mathbf{S}.  \eqno{(1.15)}\] 
  \subsection{Операторы  $S$ и  $ \mathcal{N}_d$ -  самосопряженные расширения     $\text{rot}$ и  $\nabla\text{div}$ в  $\mathbf{L}_{2}(G)$} 
  Пусть  $ \mathcal{A}_{\gamma} (G) = \{\nabla\,h, h\in H^2(G): \gamma (\mathbf{n}\cdot \nabla) h =0 \}$,  $\mathbf{A}^0_{\gamma} =\mathbf{A}^0\cap \mathcal{A}_{\gamma}$.
  
  Области определения операторов  $S$ и  $ \mathcal{N}_d$ - это пространства:
 \[   \mathbf{W}^1= \{\textbf{ f}\in  \mathbf{V}^0, \,\,  \text{rot}\textbf{ f}\in  \mathbf{V}^0 \} \quad   \text{и} \quad        \mathbf{A}^{2}=\{\textbf{ f}\in  \mathbf{A}^0_{\gamma}, \,\, \nabla  \text{div}\textbf{ f}\in  \mathbf{A}^0_{\gamma} \}, \eqno{(1.16)}\]	  	
и 
$S\mathbf{u}=\text{rot}\mathbf{u}$ при  $\mathbf{u} \in \mathcal{D}(S)=  \mathbf{W}^1,$  а
  $ \mathcal{N}_d\mathbf{v}= \nabla\text{div}\textbf{v} =  \nabla\text{div}  \nabla h$  при	
 	  	  $\mathbf{v}=\nabla h \ \in  \mathbf{A}^{2}$.
    
Согласно  оценкам  (1.12) и (1.13) при $s=0$  пространство
    $ \mathbf{W}^1\subset \mathbf{H}^{1}$ \,и  $\mathbf{A}^{2}\subset \mathbf{H}^{2}$.
    Пространство   $\mathbf{C}_0^{\infty}(G) \cap \mathbf{V}^0$ плотно в $\mathbf{V}^0$ и содержится в $ \mathbf{W}^1$;   следовательно,  $ \mathbf{W}^1$  плотно в $\mathbf{V}^0$.    
Аналогично,   $\mathbf{C}_0^{\infty} (G) \cap  \mathbf{A}^0_{\gamma}$  плотно в $ \mathbf{A}^0_{\gamma} $ и содержится в $ \mathbf{A}^{2}$;
следовательно, $\mathbf{A}^{2}$ плотно в $ \mathbf{A}^{0}$.

 Если  поля $\mathbf{u}$ и $\mathbf{v}$  в (1.14) принадлежат  $\mathcal{D}(S)$, то $\gamma  (\mathbf{n}\cdot \mathbf{u})=\gamma  (\mathbf{n}\cdot \text{rot}\,\mathbf{u})= 0$, % \newline
 $\gamma  (\mathbf{n}\cdot \mathbf{v})= \gamma  (\mathbf{n}\cdot \text{rot}\,\mathbf{v})=0$ ,  интеграл по $\Gamma$  зануляется \cite{yogi}  и  это равенсто принимает   вид:
$(\text{S} \mathbf{u}, \mathbf{v})= (\mathbf{u},\text{S}\mathbf{v})$. 

 Аналогично, если  поля  $\mathbf{u}=\nabla g$ и $\mathbf{v}=\nabla h$ в (1.15) принадлежат  $\mathcal{D}(\mathcal{N}_d)$ , то  $\gamma  (\mathbf{n}\cdot \mathbf{u})\equiv  \gamma\,(\mathbf{n}\cdot \,\nabla) g=0$,  $ \gamma  (\mathbf{n}\cdot \mathbf{v})\equiv   \gamma\,(\mathbf{n}\cdot \,\nabla) h=0$, интеграл по $\Gamma$  равен нулю и  это равенсто принимает   вид:
$
({\mathcal{N}_d} \mathbf{u}, \mathbf{v})= (\mathbf{u},{\mathcal{N}_d}\mathbf{v})$. 

Доказано, что $S$ и  $ \mathcal{N}_d$ -  самосопряженные расширения     операторов  $\text{rot}$ и  $\nabla\text{div}$ в  $\mathbf{L}_{2}(G)$
(см.  \cite{yogi,ds18}).
   \subsection{Гладкость собственных полей   операторов  $ \text{rot}$ и  $\nabla\text{div}$} 
         Спектральные задачи  для операторов  $ \text{rot}$ и  $\nabla\text{div}$ состоят в нахождении ненулевых полей $\mathbf{u}$ и $\mathbf{v}$ и чисел $\lambda$ и $\mu$  таких, что
\begin{equation*}
\klabel{srot__1_}
\text{rot} \mathbf{u}=\lambda \mathbf{u}(\mathbf{x}), \quad
\mathbf{x}\in G, \quad \gamma \mathbf{n}\cdot \mathbf{u}=0,
\quad \mathbf{u}\in \mathbf{C}^1(G)\cap \mathbf{C}(\overline{G}),
\eqno{(1.17)}
\end{equation*}
\begin{equation*}
\klabel{sgrd__1_}
\nabla\text{div}\mathbf{v}=\mu \mathbf{v}(\mathbf{x}), \quad
\mathbf{x}\in G, \quad \gamma \mathbf{n}\cdot \mathbf{v}=0,
\quad \mathbf{v}\in \mathbf{C}^2(G)\cap \mathbf{C}(\overline{G}).
\eqno{(1.18)}\end{equation*} 
Из Теорем 1,2  вытекают  важные свойства
решений  спектральных задач операторов  {\it  ротор и градиент дивергенции}: 

a)   {\it каждое ненулевое   собственное значение имеет  конечную кратность},

 b)  в  области $G$ с  гладкой границей их  обобщенные собственые поля   из $\mathbf{L}_{2}(G)
 $  является гладкими вплоть до границы. 

Доказательство.  Пусть  $\lambda\neq 0$,  а   поле $\mathbf{u}(\mathbf{x})$- соотвествующее ему решение задачи (1.17).  Такое поле $\mathbf{u}(\mathbf{x})$   есть решение однородной эллиптической задачи:
\[ \text{rot} \mathbf{u}=\lambda \mathbf{u}(\mathbf{x}), \quad \lambda\,\text{div} \mathbf{u}(\mathbf{x})=0, \quad
\mathbf{x}\in G, \quad \gamma \mathbf{n}\cdot \mathbf{u}=0,
\quad \mathbf{u}\in \mathbf{C}^1(G)\cap \mathbf{C}(\overline{G}).
\eqno{(1.19)}\]%\neq 0
Согласно Теореме 1 эта задача имеет конечное число линейно независимых решений   $\mathbf{u}_1(\mathbf{x}), ...,  \mathbf{u}_l(\mathbf{x})$, где  $l$ зависит от $\lambda$ и не зависит от  $\mathbf{u}$.  
Утверждение a) доказано. 

Решение   $\mathbf{u}(\mathbf{x})$  задачи (1.19)  принадлежит $ \mathbf{L}_{2}(G)$,  так как    $\|\mathbf{u}\|^2\equiv \int _G (\mathbf{u} \cdot \mathbf{u})d\mathbf{x} \leq V max_{\overline{G}}|\mathbf{u} \cdot \mathbf{u}|= V \|\mathbf{u} \cdot \mathbf{u}\|_{C(\overline{G})}$, где постоянная $V=\int_G 1\,d\mathbf{x}$.  % (см. \cite{vla}, Гл.IV $\S 17$).

 Согласно (1.19)  $ \text{rot} \mathbf{u}_j=\lambda \mathbf{u}_j, \,\, \text{div} \mathbf{u}_j =0$ в $ G$, %\quad
  $\gamma \mathbf{n}\cdot \mathbf{u}_j=0$.  Поэтому
 % Так как   в  $\mathbf{L}_{2}(G)$
$ \|\text{rot} \mathbf{u}_j \|= |\lambda| \| \mathbf{u}_j \|$ и
   оценка (1.12)  при   $s=0$ принимает вид: 
   $ C_0\|\mathbf{u}\|_{1} \leq (|\lambda|+1)  \|\mathbf{u}\|_{0} $, причём  постоянная $ C_0>0$. Значит,    $\mathbf{u}_j(\mathbf{x})$ принадлежит $\mathbf{H}^{1}(G)$ и
    \[ \|\mathbf{u}_j\|_{1} \leq  C_0 ^{-1}(|\lambda|+1)  \|\mathbf{u}_j
  \|_{0}, \quad    \|\mathbf{u}_j\|_0\leq \sqrt{V} \|\mathbf{u}_j\cdot \mathbf{u}_j  \|^{1/2}_{  {C}(\overline{G})}.   \eqno{(1.20)}\]
 Далее,  пусть   $s>0$ целое. Так как 
  $\| \text{rot}\mathbf{u}_j\|_s =|\lambda|\|\mathbf{u}_j(\mathbf{x})\|_s$,
 из оценки (1.12)   по индукции получаем:   
  $$\|\mathbf{u}_j\|_{s+1} \leq  C_s ^{-1} (|\lambda|+1)  \|\mathbf{u}_j\|_{s}\leq ...\leq  C_s ^{-1} ... C_0 ^{-1} (|\lambda|+1)^s  \|\mathbf{u}_j\|_{0}.             \eqno{(1.21)}$$

   Значит, поле   $\mathbf{u}_j(\mathbf{x})$ принадлежит $\mathbf{H}^{s+1}(G)$ для любого  целого  $s\geq 0$.
   
     {\it  Замечание 1.}  Известны вложения  пространств $H^{l+2}(\Omega)$  в $C^{l} (\bar {\Omega})$ при $l \geq  0$   в трехмерной области $\Omega$ и оценка 
    $ \|g\|_{C^{l} (\bar{\Omega})}\leq c_l \|g\|_{H^{l+2} (\Omega)}$ для любой функции   $g\in H^{l+2} (\Omega)$, причем постоянная 
    $c_l  > 0$ не зависит от $g$  [2, Теорема 3, § 6.2]. 
   
   Итак, поля  $\mathbf{u}_j(\mathbf{x})$ принадлежат $ \mathbf {C}^{l} (\bar G)$ для любого  целого  $l\geq 0$. Уверждение    b)  для ротора доказано.

    Аналогично, при  $\mu\neq 0$   собственное поле $\mathbf{v}(\mathbf{x})$  оператора $\nabla\text{div}$  есть решение однородной эллиптической задачи:
\[\nabla\text{div}\mathbf{v}=\mu \mathbf{v}(\mathbf{x}), \quad
\text{rot}\,\mathbf{v}=0, \quad\mathbf{x}\in G, \quad \gamma \mathbf{n}\cdot \mathbf{v}=0,
\quad \mathbf{v}\in \mathbf{C}^2(G)\cap \mathbf{C}(\overline{G}).
\eqno{(1.22)}\]
Согласно Теореме 2 эта задача имеет конечное число линейно независимых решений   $\mathbf{v}_1(\mathbf{x}), ...,  \mathbf{v}_k(\mathbf{x})$, где  $k$ зависит от $\mu$ и не зависит от  $\mathbf{v}$.  Утверждение a) доказано. \quad   
Любое решение   $\mathbf{v}_j(\mathbf{x})$  задачи (1.22)  принадлежит $\mathbf{L}_{2}(G)$,  так как   $\|\mathbf{v}\|^2_{\mathbf{L}_{2}(G)}\leq V\|\mathbf{v} \cdot \mathbf{v} \|_{  {C}(\overline{G})}$, где постоянная $V=\int_G 1\,d\mathbf{x}$.

 Ввиду того, что
   $ \|\nabla\text{div}\mathbf{v}\|= |\mu| \|\mathbf{v}\|$   в  $\mathbf{L}_{2}(G)$,
оценка (1.13)  при   $s=0$ принимает вид: 
$ C_0\|\mathbf{v}\|_{2} \leq (|\mu|+1)  \|\mathbf{v}\|_{0} $, причём  постоянная $ C_0>0$.

 Значит,    $\mathbf{v}_j(\mathbf{x})$ принадлежит $\mathbf{H}^{2}(G)$, и
\[ \|\mathbf{v}_j\|_{2} \leq  C_0 ^{-1}(|\mu|+1)  \|\mathbf{v}_j
\|_{0}, \quad    \|\mathbf{v}_j\|_0\leq \sqrt{V} \|\mathbf{v}_j\cdot \mathbf{v}_j  \|^{1/2}_{  {C}(\overline{G})}.    \eqno{(1.23)}\]
Далее,  пусть   $s>0$ целое. Так как 
$\|\nabla\text{div}\mathbf{v}_j\|_s =|\mu|\|\mathbf{v}_j(\mathbf{x})\|_s$,
из оценки (1.13)   по индукции получаем:   
$$\|\mathbf{v}\|_{2s+2} \leq  C_{2s} ^{-1} (|\mu|+1)  \|\mathbf{u}\|_{2s}\leq ...\leq  C_{2s} ^{-1} ... C_0 ^{-1} (|\mu|+1)^s  \|\mathbf{u}\|_{0}.      \eqno{(1.24)} $$

Значит,    $\mathbf{v}_j(\mathbf{x})$ принадлежит $\mathbf{H}^{2s+2}(G)\subset \mathbf{C}^{2s} (\bar {G}) $ для любого  целого  $s\geq 0$.
Утверждение b) доказано.

 \subsection{Гладкость базисных  полей   пространств $\mathcal{A}_H$ и $\mathcal{B}_H$}

Пространства   $\mathcal{A}_H$ и  $\mathcal{B}_H$  определяются решениями эллиптических систем (1.4)  и  (1.5)  в $ \mathbf{L}_{2}(G)$. Из формул (1.8) видно, что  компоненты  этих решений являются гармоническими функциями,  а значит,  они имеют непрерывные производные любого порядка.  Это впервые заметил  Герман Вейль для решений системы  (1.5) (\cite{hw} Теорема 1).
Краевые задачи   (1.4)  и  (1.5) удовлетворяют  условиям    В.Солонникова в Теореме 1.1  работы \cite{so71}. Откуда
получаем, что  пространства   $\mathcal{A}_H$ и  $\mathcal{B}_H$  конечномерны и их базисные поля  $\mathbf{g}_i(\mathbf{x})$ и $\mathbf{h}_j(\mathbf{x})\in 
\mathbf{C}^{\infty} (\bar G)$,  $i=1,..,\rho_1<\infty $, $j=1,..,\rho<\infty$.    
Отметим, что для $\mathbf{g}_i(\mathbf{x})$ и $\mathbf{h}_j(\mathbf{x})$ имеются оценки вида  (1.21) с   $\lambda=0 $.   
Borchers W.,  Sohr  H.  доказали \cite{boso}, что число  $\rho$ есть род границы $\Gamma$ области  $G$. 
В частности, если область  $\Omega$  гомеоморфна шару,  то $\rho=0$.

%{\it Замечание. } Наша гипотеза состоит в том, что  $\rho_1 = 0$ в %области  $\Omega$,  гомеоморфной шару.  Пока, это доказано только %для шара.

Если область  $\Omega$  гомеоморфна шару,  а $\mathbf{u}$  -решение задачи  $(1.5)$,    определяющей   $\mathcal{B}_H$, то   $\mathbf{u}= \nabla\, h$,  а функция $h$ - решение задачи Неймана для оператора Лапласа: $\Delta\,h =0 $ в $\Omega$,  $ \gamma (\mathbf{n}\cdot \nabla)\, h=0$.
Решение этой задачи   $N$ есть произвольная постоянная  $h=Const$. 
Следовательно,   $\mathbf{u}\equiv 0$ и пространство  $\mathcal{B}_H$ пусто.

%Определим пространство  $\mathcal{A}_H$ 
Решение задачи (1.25)  в шаре $B, |\mathbf{x}|<R$,  сводится к задаче: $\Delta\,h =C$ в $B$,  где $C$ -произвольная постоянная,  с условием  Неймана      $\gamma (\mathbf{n}\cdot \nabla)\, h=0$.  Пусть $\mathbf{r}=\mathbf{x}$ - радиус-вектор, тогда нормаль  $\mathbf{n}=\mathbf{x}/r$   на границе шара, где
  $r= |\mathbf{x}|=R$.  Частное решение уравнения Пуассона  $\Delta\,h =C$ имеет вид: $h= 1/6C |\mathbf{x}|^2= 1/6Cr^2 $.
  Дифференцируя  по $r$,  получаем  $\gamma (\mathbf{r}\cdot \nabla)\, h= 1/3Cr|_{r=R}= CR/3$. Граничное условие Неймана принимает вид: $ CR/3=0$.  Значит,  $ C=0$ и  пространство 
 $\mathcal{A}_H(B)$ в шаре $B$ пусто.

 \subsection{Ортогональные базисы  в  $\mathcal{A}$, \, $\mathcal{B}$\,  и   в  $\mathbf{L}_{2}(G)$} Пространство $\mathbf{A}^2$ плотно в $\mathbf{A}^0 $ и
$\mathbf{A}^2\subset \mathbf{H}^2$.
Собственные поля $\mathbf{q}_{j}(\mathbf{x})$  оператора  $\nabla\mathrm{div}$ с ненулевыми собственными значениями ${\mu}_{j}$  принадлежат $\mathbf{A}^2$.

 {\small Множество собственных значений $\mu=-\nu^2$ этого оператора 	счётно, отрицательно и каждое из них имеет
	конечную кратность. Перенумеруем их в порядке возрастания их модуля: 	$0<-\mu_1\leq -\mu_2\leq ...$,  повторяя $ \mu_k$ столько раз, какова	его кратность. Соотвествующие вектор-функции  обозначим через
	$\mathbf{v}_{1}, \mathbf{v}_{2}$, ..., так чтобы каждому
	значению $\mu_{k}=-\nu^2_k$ соответствовала только
	одна  функция $\mathbf{v}_{k}$: $\nabla \mathrm{div}
	\mathbf{v}_{k}=-\nu^2_k \mathbf{v}_{k}$,  $ \gamma\mathbf{n}\cdot\mathbf{v}_{k}=0, $
	$k=1,2,...$.	
	Собственные функции, соответствующие одному и тому же
	собственному значению, выберем ортонормальными, используя процесс
	ортогонализации Шмидта  (см. \cite{vla}). 
	Поля, соответствующие
	различным с.- значениям, ортогональны. Их нормируем.
	Нормированные собственные поля  градента дивергенции обозначим 
	$\mathbf{q}_{l}$,   \quad  $ l=1,2,...$,  норма  $\|\mathbf{q}_{l}\|=1$.
	Они составляют полный ортонормированный базис в классе 
$\mathbf{A}^{0}$. 
  Зафиксируем его.}

Аналогично строится базис в классе $\mathbf{V}^{0}$ 
 \cite{saVS20}.
 
 {\it Замечание.}  Согласно (1.9)  оператор $\Delta \mathbf{u} \equiv -\mathrm{rot}^2 \mathbf{u}$  при  $\mathbf{u}\in \mathcal{B}$.   
 Собственные векторы ротора    всегда  встречаются парами: 
  каждому с.-полю    $\mathbf{u}^{+}_{j}$ с  $\lambda_j>0$  соответствует с.-поле   $\mathbf{u}^{-}_{j}$ с  $-\lambda_j$.  	
Это их свойство  в \cite{yogi} не отмечено.

Зафиксируем в  $\mathbf{V}^{0}$ ортонормированный базис $\{\mathbf{q}^{+}_{j}, \mathbf{q}_{j}^{-}\},  \quad  \|\mathbf{q}^{\pm}_{j}\|=1$: 
\[\mathrm{rot}
\mathbf{q}_{j}^{\pm}=\pm\lambda_j\, \mathbf{q}_{j}^{\pm}, \quad   \gamma\mathbf{n}\cdot\mathbf{q}_{j}^{\pm}=0,   \quad   j=1, 2, ...,  \quad  \mathbf{q}^{\pm}_{j}\in \mathbf{C}^\infty(\bar{G}).     \eqno{(1.26)} \]
Учитывая базисы пространств $\mathcal{A}_H$, \, $\mathcal{B}_H$
согласно (1.3), (1.6) видим, что  объединение $\{g_l\}$,	$\{\mathbf{q}_{l}\}$,  $\{h_j\}$  и $\{\mathbf{q}^{+}_{j}, \mathbf{q}_{j}^{-}\}$ есть базис объемлющего пространства  $\mathbf{L}_{2}(G)$.

  Спектральные задачи для операторов      	ротор и градиент дивергенции  в шаре    решены автором   полностью  в   \cite{saUMJ13}. 
   
     \subsection{Спектральная задача для оператора ротор в шаре  $B$.  Решение}  	
   
    Имеется несколько способов   решения  этой задачи  \cite{wo58, cdtgt,saUMJ13}.
   
   Учитывая  приложения \cite{wo58} и конкурирующие интересы           \cite{ cdtgt},  кратко изложим наш путь  решения   этой задачи  \cite{saUMJ13}.

  	Собственные числа $\lambda_{n,m}$  ротора  в шаре радиуса $R$ равны $\pm \rho_{n,m}/R$,  где числа $\pm \rho_{n,m}$ - нули функций  $\psi_n(r)$:
  	\begin{equation*}\label{bes  21}
  	\psi_n(z)
  	= (-z)^n\left(\frac d{zdz}\right)^n\left(\frac{\sin z}z\right),  \quad  m,n\in {\mathbb {N}}. \eqno{(1.27)}
  	\end{equation*}
  Функции  $\psi_n(z)$- это цилиндрические  функции $J_{n+1/2}(z)$,  где  $n\geq 0$ целое.  Это заметил ещё Леорнард Эйлер (см. \cite{vla} $\S23$).
  Числа $\pm \rho_{n,m}$ и  $ \rho^2_{n,m}>0$ - нули функций  $\psi_n(z)$.
  Кратность  собственного значения  $\lambda^{\pm}_{n,m}$ равна $2n+1$. 
  
  Пусть $\mathbf{i}_r, \mathbf{i}_\theta,
  \mathbf{i}_\varphi$-репер, поле  $\mathbf{u}=u_r\,\mathbf{i}_r+ u_{\theta} \mathbf{i}_{\theta}+ u_{\varphi}
  \mathbf{i}_\varphi$.% \newline  

  	{\bf Формулы решений задачи (1.26).}
    		 {\it  Ненормированные собственные поля
    		 	 $ \mathbf{u}_{\kappa
  		}^{\pm }$ задачи  (1.26) в сферических координатах вычисляются по
  		формулам:}
  		\begin{equation*}
  		\label{vrsff   2}\begin{array}{c}
   \mathbf{u}_{\kappa}^{\pm }=c_{\kappa }^{\pm }(\pm\lambda _{n,m}
  		r)^{-1}{\psi }_{n}
  		(\pm\lambda _{n,m}r)Y_{n}^{k}(\theta ,\varphi  )\,\mathbf{i}_r+\\
  		c_{\kappa }^{\pm }{{(\pm\lambda _{n,m}^{\pm }r)}^{-1}}
  		Re[\Phi _{n}(\pm\lambda _{n,m}r)](Re HY_{n}^{k}\,
  		\mathbf{i}_\varphi+
  		Im HY_{n}^{k}\,\mathbf{i}_\theta)+\\
  		c_{\kappa }^{\pm }{{(\pm\lambda _{n,m}r)}^{-1}}
  		Im[\Phi _{n}(\pm\lambda _{n,m}r)](-Im HY_{n}^{k}\,
  		\mathbf{i}_\varphi+
  		Re HY_{n}^{k}\,\mathbf{i}_\theta).\end{array}   \eqno{(1.28)}\end{equation*}
  		где
  		$Y_{n}^{k}(\theta ,\varphi  )$--сферические функции, числа
  	$c_{\kappa }^{\pm }\in \mathbb{R}$ -произвольны, $\kappa=(n,m,k)$- мульти-индекс,
  $m{{,}^{{}}}n\in \mathbb{N}$,    $|k|\le n$, 
   	\[\Phi _{n}(\lambda\,r)=
  	\overset{\mathop{{}}}\,\int\limits_{0}^{r}{}\,{{e}^{i\lambda   			(r-t)}}^{{}}
  	{\psi }_{n}(\lambda t){t}^{-1}dt,
  	\quad Im \Phi _{n}(\pm\rho_{n,m})=0,        \eqno{(1.29)}\]
  	\begin{equation*}
  	\label{oph   1}{2k}
  	\text{H}v=
  	\left(\sin ^{-1}\theta {\partial }_{\varphi }+
  			i{\partial }_{\theta } \right)v,  \quad 	\text{K}w=
  		\sin ^{-1}\theta 	\left({\partial }_{\theta } \sin \theta+
  			i{\partial }_{\varphi } \right)w .   \eqno{(1.30)} \end{equation*}  
  Решению  	этой	   спектральной задачи
 способствовали  наблюдения автора:

    1. {\it 	Пусть поле $\mathbf{u}$ -  решение  спектральной задачи (1.19) в шаре $B$,  	$\mathbf{x}$-радиус- вектор, а
  	  $v(\mathbf{x})$- их  скалярное  		произведение   $\mathbf{x}\cdot \mathbf{u}=r\,u_r$.
        Тогда	функция    $v(\mathbf{x})$  есть решение  спектральной задачи Дирихле для уравнения Лапласа: }
  	\[	\label{ldo__1_}           -\Delta v=\lambda^{2}\,v  \quad \text{в} \quad B,
  \quad v|_{S}=0,   \quad \text{с условием} \quad  v(0)=0.                 \eqno{(1.31)}\]

  		2. {\it  Уравнения $\mathrm{rot}\mathbf{u}=\lambda \mathbf{u}, \,  \mathrm{div}\mathbf{u}=0$, записанные   в сферических координатах,  представляются 	в виде двух комплексных уравнений   			
  			\[({\partial}_r - i\lambda) r\,w= r^{-1}H\,v,
  			\quad  K\,w=\lambda\,v-i\,r^{-1\,}\,{\partial}_r (r\,v),    \eqno{(1.32)} \]
  			относительно функций $v=ru_r$ и $w=u_{\varphi}+iu_{\theta}$	 с операторами  $H$  и $K$  в (1.30).  
  		 }
  	 
  	 3.  {\it   Уравнения (1.31) на  функцию $v$ являются условиями   	 совместности системы (1.32).}
  
  		Таким образом,  решение задачи сводится к решению:
  		  		
  	{ \small $1^0)$  спектральной задачи Дирихле - Лапласа  (1.31).  Её решения- пары  $ \lambda_{\kappa}^2=(\rho_{n,m}/R)^2$
  		и  ${v}_{\kappa}=c_{\kappa }{\psi }_{n}(\rho_{n,m}r/R )Y_{n}^{k}(\theta ,\varphi  )$, такие что ${\psi }_{n} (\rho^2 _{n,m})=0$
  		. Условие 	  v(0)=0    обеспечивается обнулением постоянных $c_{\kappa }=0$  при $\kappa=(0,m,0)$ (см. В.С.Владимилов   \cite{vla} гл.V $\S26$).  Они определяют  $ \lambda^{\pm}_{\kappa}=\pm \rho_{n,m}/R$-собственные  значения задачи (1.17) в шаре $B$  и  функции    $ u_{r,\kappa }=v_{\kappa}/r $-  радиальные компоненты собственных полей. }
  	
  	{ \small $2^0)$ к интегрироваию уравнений (1.32) с  $ \lambda= \lambda^+_{\kappa}>0$ и $v=v^+_{\kappa}$, а затем с  $ \lambda= \lambda^-_{\kappa}<0$ и $v=v^-_{\kappa}$,  и  вычилению комплексных функций $w^{\pm}_{\kappa}$,     задающих  касательные компоненты полей   $\mathbf{u}^{\pm}_{\kappa}$; они 
  		определятся однозначно  условием: $w^{\pm}_{\kappa}\in L_2(B)$
  		
  		\small $3^0)$ к 	построению полей 	$\mathbf{u}_{\kappa}^{\pm}(\mathbf{x})\in \mathbf{L}_2(B)$. }    	
  	
  	В итоге, получаем список решений (1.28).
  	
 {\it Замечание.}  	Позже 	уравнение (1.31) на функцию  $v=r\,u_r$  при  минимальном  собственном значении $\lambda=4.4934.../R$ автор обнаружил  в статье \cite{cdtgt}.

   \subsection{Спектральная задача для оператора   $\nabla\mathrm{div}$ в шаре  $B$.  Решение}      Собственные значения  оператора $-\nabla\mathrm{div}$ равны $\nu_{n,m}^2$, где $\nu_{n,m}=\alpha_{n,m}/R$,  а  числа  $\alpha_{n,m}$ - нули производных  $\psi'_n(r)$,   $n \geq 0,\, m\in {\mathbb {N}}$ ;     	
 	кратность  собственных значений  $\nu^2_{n,m}$    равна $2n+1$.	 	
 
 	Собственные поля $\mathbf{v}_{\kappa}= \nabla  g_{\kappa}$    градиента дивергенции  - решения задачи:
 	\[  -\nabla \mathrm{div}
 	\mathbf{v}_{k}=\nu^2_{\kappa}\mathbf{v}_{\kappa}, \quad  \gamma\mathbf{n}\cdot\mathbf{v}_{\kappa}=0,  \quad  \mathbf{v}_{\kappa}= \nabla  g_{\kappa}\in \mathcal{C}^\infty(\bar{G}).   \eqno{(1.33)} \] 
Так как   $ \nabla \mathrm{div}\nabla  g_{\kappa}\equiv  \nabla 
\Delta_c\,g_{\kappa}= 
\Delta_c \, (\nabla g_{\kappa}) =  -\nu^2_k \,(\nabla g_{\kappa}),  \quad \gamma(\mathbf{n}\cdot  \nabla ) g_{\kappa}=0$ эта задача сводится к задаче Неймана для скалярного оператора Лапласа и градиенту фукций $g_{\kappa}$.  Матричный (3 x 1) оператор   $ \nabla \mathrm{div}\nabla  g \equiv  \nabla 
\Delta_c\,g$ эллиптичен.

Соответствующие  $\nu^2_{\kappa}\equiv  \nu^2_{n,m}$ собственные функции  $g_{\kappa}$ имеют вид:
\[ g_{\kappa }(r, \theta ,\varphi)=c_{\kappa }\psi _{n}
(\alpha _{n,m}r/R)Y_{n}^{k}(\theta ,\varphi  )    \eqno{(1.34)} \]

Поля   $\mathbf{v}_{\kappa}= \nabla  g_{\kappa}$  являются решениями задачи (1.25); их компоненты $(v_r, v_{\theta}, v_{\varphi })$  имеют вид:
 $ v_{r, \kappa }(r, \theta ,\varphi)=c_{\kappa } (\alpha _{n,m}/R)\psi'_{n} (\alpha _{n,m}r/R)Y_{n}^{k}(\theta ,\varphi), $ 
 \[ (v_{\varphi}+iv_{\theta})_{ \kappa}=c_{\kappa } (1/r)       \psi_{n} (\alpha _{n,m}r/R)\,H\,Y_{n}^{k}(\theta ,\varphi  )  \eqno{(1.35)}                  \]
 
 При $ \kappa=(0,m,0)$ функция $Y_{0}^{0}(\theta ,\varphi)=1, HY_{0}^{0}(\theta ,\varphi)=0,$  поэтому 
 
\[ v_{r, (0,m,0) }(r)=  c_{ (0,m,0)}(\alpha _{0, m}/R) \psi'_{0} (\alpha _{0,m}r/R) , \quad  (v_{\varphi}+iv_{\theta})_{ (0, m, 0)}=0.       \eqno{(1.36)}             \]

  Построенный в  шаре  $B$  базис из собственных полей операторов градиента дивергенции и  ротора является полным в $\mathbf{L}_{2}(B)$,
		так как $\mathbf{L}_{2}(B)=\mathbf{A}^0 \oplus \mathbf{V}^0$.

 \subsection{Потоки с  минимальпой энергией,  изображения}
 	Эти формулы используются при рассчетах поля скоростей $\mathbf {u}_{\kappa }^{\pm}(\mathbf {x})$ и визуализации вихревых потоков.  %при заданном $\kappa$. 
 Формулы  полей  $\mathbf{u}_{\kappa}^{\pm}(\mathbf{x})$  при  $n=1$, $\kappa=(1,1,0)$ и  $\kappa=(1,1,\pm 1)$ выражаются  наиболее просто. Так,  компоненты поля $\mathbf{u}^{+}_{(1,1,0)}(\mathbf{x})$ имеют вид:    % \eqno{(1.28)}
 \[u_r=2\rho( r\rho)^{-3}( sin\,(r\rho)- r\rho cos\,(r\rho)) cos\,\theta,   \] 
 \[u_\theta =( r\rho)^{-3}( sin\,(r\rho)- r\rho \,cos\,(r\rho)-( r\rho)^{2}sin\,(r\rho)) sin\, \theta ,       \eqno{(1.37)} \]
 \[u_{\varphi} = ( r\rho)^{-2}(( sin\,(r\rho)- r\rho \,cos\,(r\rho)) sin\, \theta.  \]  
 	Профессор Исламов Г.Г.\cite{gais}, используя эти формулы   и  программу  Wolfram Mathematica  осуществил
 визуализацию линий тока поля   $u_{1,1,0}^{+}(\mathbf {x})$ ротора  радиуса 1 со значением $\rho =\rho_{1,1}$ \footnote{ 
 	http://www.wolfram.com/events/ technology-conf.-ru/  2016/resources.html }.   
 Траектория движения  трёх соседних точек напоминает ленту, которая 	огибает тороидальную катушку как волна  (см. также  катушка Исламова  в \cite{saVS20}).
  \footnote{Исламов Галимзян Газизович (02.02.1948-22.11.2017) }
  
 	В связи с задачами  астрофизики  S. Chandrasekhar, P.С. Kendall изучали  собственные поля ротора в шаре \cite{chake}  и  в цилиндре. 
 	 Они нашли   элементарный способ их вычисления  в цилиндре  (с условием периодичности вдоль оси).
 	
 	D.Montgomery, L.Turner, G.Vahala,  	 изучая магнито гидродинамическую турбулентность в цилиндре  \cite{motva},
 	 исползовали эти формулы.
 	 
 	  Они показали, что три интегральных инварианта имеют простые квадратичные выражения в тестатье \cite{cdtgt}.
 Они приводят  также   соответствующую	 $\lambda_{1,1}$  формулу  собственного поля ротора в шаре  (см. Theorem A). К сожалению,  с опечаткой... ( Может это код?). 
 %$1/\lambda$  вместо $\lambda$.  
 
 Исправив её,   мы пришли к    компонентам   (1.37)  поля  $\mathbf{u}^{+}_{(1,1,0)}(\mathbf{x})$. 

В Fig.1  в  \cite{cdtgt}  представлены 
интегральные кривые  поля $\mathbf{u}^{+}_{(1,1,0)}(\mathbf{x})$    и даётся их описание. \quad  Цитирую: "они заполняют семейство концентрированных "торов"  с  замкнутой орбитой  "ядра,"  \, типичных  для осесимметричных собственных полей ротора; специальная орбита  начинается на южном полюсе сферы в момент  времени $-\infty $,  проходит вертикально вверх по оси z и достигает северного полюса ко времени  $+\infty$; орбиты на граничной сфере   начинаются на северном полюсе в момент времени  $-\infty $,  продолжаются по линиям долготы к южному полюсу до момента времени  $+\infty $; имеются две стационарные точки в её полюсах."

Мы провели независимое исследование этого поля. Галимзян Исламов демонстрировал эти орбиты  вживую на  экране в  МГУ  во время нашего совместного доклада на конференции  Бицадзе-100,  факультет ВМК 2016.

 Авторы  \cite{cdtgt}  приводят также формулы   базисных  полей ротора   в шаре  для других собственных значений.  
 Для сравнения  мы  приводим  весь список  (1.28).

   \subsection{Степени оператора Лапласа в классах $\mathcal{A}$, \, $\mathcal{B}$\,  и   в  $\mathbf{L}_{2}(G)$}
Из формул (1.9) при $k=2, 3, ...$ имеем
\[\Delta^k \mathbf{v}\equiv  (\nabla \mathrm{div})^k \mathbf{v}   \quad  \text{при}    \quad  \mathbf{v}\in \mathcal{A},  \quad  \Delta^k \mathbf{u} \equiv (-1)^k(\mathrm{rot})^{2k}\,\mathbf{u}  \quad  \text{при}  \quad      \mathbf{u}\in \mathcal{B}.   \eqno(1.38)\]
В  $\mathbf{L}_{2}(G)$ оператор $\Delta^k$
выражается через  $(\nabla \mathrm{div})^k$ и $(\mathrm{rot})^{2k}$,   а также через скалярный оператор $\Delta_c^k = (\partial_1^2+\partial_2^2+\partial_3^2)^k$ : 
\begin{equation*}\klabel{dd 1}\mathrm{\Delta}^k\,\mathbf {v}
= (\nabla \mathrm{div})^k\,\mathbf {v} + (-1)^k\,(\mathrm{rot})^{2k}\, \mathbf {v}= \Delta^k_c\, I_3\,\mathbf {v} , \quad \text{где}
\,\,  \mathbf{v}=(v_1, v_2, v_3). \eqno(1.39) \end{equation*} \quad  
Эти формулы следуют из формул (1.8),  учитывая, что операторы $\mathrm{rot}$ и $\nabla \mathrm{div}$ аннулируют друг друга.  Они являются проекторами :
 $\nabla\mathrm div$ проектирует $\mathbf{L}_{2}(G)$ на $\mathcal{A}$, а  $\mathrm{rot}$ -   на $\mathcal{B}$   .

Как мы отечали во введении,  в  периодическом случае С.Л.Соболев  доказал      Т е о р е м у  XII.13.  Мы рассмотрим её аналоги.

Операторы   $(\nabla \mathrm {div})^p$ и $(\mathrm {rot})^{2q}$, где  $ p$ и  $ q$ - натуральные  числа,
- аналоги  полигармонических  операторов $\Delta^m$
 в классах $\mathcal{A}$   и  $\mathcal{B}$ (см.  формулы  (1.38)).
Мы покажем,  что 
 оператор $(\nabla \mathrm {div})^{2p}$ переводит
произвольное поле $\mathbf{w}$  из   $ {A}^{2p}$ в 
$(\nabla \mathrm {div})^{2p}\,\mathbf{w}=\rho$ --  элемент  ${A}^{-2p}\equiv  ({A}_0^{2p})^*$  (см. п. 2. 3);

а оператор $(\mathrm {rot})^{2q}$ переводит
 произвольное поле $\mathbf{u}$  из   $ {W}^{q}$ в 
$(\mathrm {rot})^{2q}\,\mathbf{u}=\mathbf{v}$ --  элемент  ${W}^{-q}\equiv  ({W}_0^{q})^*$ (см. п.3.6).

Доказаны и обратные утверждения (см. Теоремы 4 и 5).

  \section{Пространство $\mathcal{A}$ потенциальных полей }

 В статье автора \cite{ds18} детально рассмотрена структура класса $\mathcal{A}$ потенцциальных полей, его базис    и оператор $ \mathcal{N}_d$.  Здесь мы рассмотрим  его подпространства  $\mathbf{A}^{2k}$.
    По  определению ${\mathcal{{A}}}(G) =\{\nabla h, h\in H^1\}$,  $\mathcal{A}_H $ - ядро оператора  $ \nabla\mathrm{div}$  в $\mathcal{A}$,  а    $\mathbf{A}^0$ - его ортогональное дополнение в  $\mathcal{A}$,  $\mathcal{A}=\mathcal{A}_H\oplus \mathbf{A}^0$,    \quad
     $ \mathcal{A}_{\gamma} (G) =\{\nabla h, h\in H^2(G): \gamma(\mathbf{n}\cdot \nabla) h=0\}$,   \quad
    $ \mathcal{A}^0_{\gamma}  = \mathbf{A}^0\cap \mathcal{A}_{\gamma}$
   
  Подпространство
$\mathbf{A}^2=\{ \mathbf{v}\in
\mathcal{A}^0_{\gamma}: \nabla\mathrm{div} \mathbf{v}\in
\mathcal{A}^0_{\gamma}\}$  есть область определегия оператора   $ \mathcal{N}_d$;  оно плотно в $\mathbf{A}^0$ и
$\mathbf{A}^2\subset \mathbf{H}^2$ (согласно п. 1.5). 
Собственные поля $\mathbf{q}_{j}(\mathbf{x})$  оператора  $\nabla\mathrm{div}$ с ненулевыми собственными значениями $(-\nu^2_j)$:
 $\nabla \text{div}\, \mathbf{q}_{j}= -\nu^2_j\,\mathbf{q}_{j}$, \,
$\gamma(\mathbf{n}\cdot \mathbf{q}_{j})=0$,
принадлежат пространству  
$\mathbf{A}^2$.
Они составляют ортонормальный базис $\{\mathbf{q}_{j}(\mathbf{x})\}$  в 	$\mathbf{A}^{0}$ 
 (см. п. 1.8).  
  	Проекция поля $\mathbf{f}\in \mathbf{L}_2(G)$ на
$\mathbf{A}^0$ имеет вид:  \[\mathcal{P}_{\mathcal{A}}\mathbf{f}\equiv 	\mathbf{f}_{\mathcal{A}}(\mathbf{x})=	
\lim_{n\rightarrow\infty}(\mathbf{f}^n_{\mathcal{A}})=\sum_{j=1}^\infty
(\mathbf{f},\mathbf{q}_{j})\mathbf{q}_{j}(\mathbf{x}),   \eqno{(2.1)}\]	
где	 $\mathbf{f}^n_{\mathcal{A}}$--  частичные суммы  этого
ряда. 

Оператор $ \mathcal{N}_d$  определен %на $ \mathbf{A}^{2}$ 
  и совпадает с $\nabla\mathrm{div}$ на  $ \mathbf{A}^{2}$, поэтому	
\[ \mathcal{N}_d\mathbf{f}_{\mathcal{A}}=
\lim_{n\rightarrow\infty}
\nabla\mathrm{div}\,(\mathbf{f}^n_{\mathcal{A}})
=-\sum_{j=1}^{\infty}\nu^2_j
(\mathbf{f},\mathbf{q}^{}_{j})\mathbf{q}^{}_{j}(\mathbf{x}),   \eqno{(2.2)}\]
если ряд	сходится и принадлежит $\mathbf{A}^0$.	
Это так, если $f\in\mathbf{H}^2(G)$.  

Доказано, что 		оператор $\mathcal{N}_d$ замкнут и самосопряжён  \cite{ds18}.  

 \subsection{Подпространства  $\mathbf{A}^{2k}$ в $\mathcal{A}$}
Рассмотрим ещё пространства\footnote{Они совпадают с пространствами $\mathbf{A}^{2k}_ {\gamma}$  в  \cite{ds18}, если пространство  $\mathcal{A}_{H}$ пусто. }
\[\mathbf{A}^{2k}=\{\textbf{ f}\in  \mathcal{A}^0_ {\gamma},..., (\nabla  \text{div})^k\textbf{ f}\in \mathcal{A}^0_ {\gamma} \}, \quad k=1,2, ... \eqno{(2.3)}\]	  		  	  		  	 
{\it Замечание. }  Согласно  оценке (1.13) пространство $\mathbf{A}^{2k}\subset \mathbf{H}^{2k}$. 
Оно    является   проекцией пространства Соболева  $\mathbf{H}^{2k}$ порядка ${2k} $ на класс  $\mathcal{A}$, так как  для любого поля   $\mathbf{f}\in  \mathbf{H}^{2k} $  его проекция $\mathcal{P}_{A}\mathbf{f}\in  \mathbf{A}^{2k}$;
 если  же $\mathbf{f}\in   \mathbf{A}^{2k}$,   то $\mathcal{P}_{A}\mathbf{f}= \mathbf{f}$,  а 	его проекция 
  на 	  	  $\mathcal{B}$ равна 0.  

	  Пространство $\mathcal{A}^0_{\gamma}$ ортогонально ядру оператора
$\mathcal{N}_d$
в $\mathbf{L}_{2}(G)$, поэтому   $\mathcal{N}_d$ имеет единственный обратный	  	  оператор: 
\begin{equation*}\label{Nobr__2_}\mathcal{N}_d^{-1}\mathbf{f}_{\mathbf{A}}=
-\sum_{j=1}^{\infty}\nu_j^{-2}
(\mathbf{f},\mathbf{q}_{j})\mathbf{q}_{j}(\mathbf{x}).    \eqno{(2.4)}\end{equation*} 
	Оператор $\mathcal{N}_d^{-1}$ - компактен. 

Следствие. {\it Спектр оператора $\mathcal{N}_d^{-1}$ точечный с
	единственной точкой  накопления в нуле,\quad
	$\nu^{-2}_j\rightarrow 0$ при ${j\rightarrow\infty}$.}

\subsection{Сопряжённые пространства  $\mathbf{A}^{-2k}$}
По определению пространство %   \quad 
 ${H}^{s}_0(G)$ есть замыкание	в норме ${H}^{s}(G)$  функций из ${C}^{\infty}_0(G)$.  
    $\mathcal{A}_0 =\{\nabla h, h\in H^1_0\}$, \quad
$\mathbf{A}^{2k}_ {0}=\{\textbf{ f}\in  \mathcal{A}_ {0}, ..., (\nabla   \text{div})^k\textbf{ f}\in \mathcal{A}_ {0} \}$.  \quad
 Пространство линейных непрерывных функционалов над   $\mathbf{A}^{2k}_ {0}$,  обозначим $(\mathbf{A}^{2k}_0)^*$. Они равны нулю на $ \mathcal{A}_H$ (см. п.2.3).
 
 В п. 2.4 мы покажем, что эти пространства можно отождествить с пространствами $\mathbf{A}^{-2k}$ 
  порядка $-2k$.
   Наконец, $ \mathcal{A}^*$- это объединение  $\mathbf{A}^{-2k}$ при  $k\geq 1$.
 
Цепь вложений пространств $\mathbf{A}^{2k}$ 
имеет вид: \[\subset \mathbf{A}^{2k}\subset...\subset\mathbf{A}^{2}\subset  \mathbf{A}^0\subset \mathbf{A}^{-2}\subset...\subset \mathbf{A}^{-2k}\subset \eqno{(2.5)}\]

Операторы $\mathcal{N}_d:\mathbf{A}^{2k} \rightarrow \mathbf{A}^{2(k-1)}$ \, обратимы  при    $k> 1$   и 
$$\|\mathcal{N}_d^{-1}\mathbf{f}\|^2_{\mathbf{A}^{2k}} 
\leq  c^2_k \|\mathbf{f}\|^2_{\mathbf{A}^{{2(k-1)}}},  \quad
\|\mathcal{N}_d  \mathbf{f}\|^2_{\mathbf{A}^{2(k-1)}}\leq  
c^{-2}_k \|\mathbf{f}\|^2_{\mathbf{A}^{2k}},       \eqno{(2.6)}$$
  где
$c^2_k= max _j(1+1/{\nu}_{j}^{2k})$,   а $1/{\nu}_{j}\to 0$ при $j\to \infty $. 

{\it Замечание.} Автор изучал также оператор  $ \mathcal{N}_d+\lambda I$ в  \cite{saVS20, ds18}, доказана

\begin{theorem} \klabel{Nd_3}Оператор $\mathcal{N}_d+\lambda I:\mathbf{A}^{{2(k+1)}} \rightarrow \mathbf{A}^{2k}$ -фредгольмов при $k\geq0$. 
	Если $\lambda \overline{\in} Sp (\mathcal{N}_d)$,  ,
	то 	 оператор  $ \mathcal{N}_d+\lambda I$  (и его обратный)  отображает    пространство  $ \mathbf{A}^{2(k+1)}$
	на  $  \mathbf{A}^{2k}$ (и обратно) 	взаимно	однозначно и непрерывно.\end{theorem}

Оператор     $ \mathcal{N}_d \mathbf{u}$ совпадает с $\nabla \mathrm{div}\mathbf{u}$, если   $\mathbf{u}\in \mathbf{A}^{2}\equiv \mathcal{D}(\mathcal{N}_d) $. Поэтому оператор $(\nabla \mathrm{div})^{k}$ на $\mathbf{A}^{2k}\subset \mathbf{A}^{2}$ совпадает с $ \mathcal{N}_d^{k}$ при    $k> 1$. 
 
%\end{theorem} \begin{theorem}
\subsection{Оператор    $\mathcal{N}^{2k}_d$ в пространстве $\mathbf{A}^{2k}$}Основное утверждение.

  {\it Оператор    $\mathcal{N}^{2k}_d$ отображает пространство $\mathbf{A}^{2k}$  на  $\mathbf{A}^{-2k}$ и обратно.}

Этапы доказательства:

 Шаг 1-й:  {\it Оператор $ \mathcal{N}_d^{2k}$ отображает пространство
$ \mathbf{A}^{2k}$  на   $ (\mathbf{A}^{2k}_0)^*$.}

Действительно,  	пусть $\mathbf{w}$ произвольный элемент из $\mathbf{A}^{2k}$, а
$\mathbf{w}_{\eta}$ -- средняя вектор-функция для него,
$\mathbf{w}_{\eta}\in\mathbf{A}^{2k}_0$;   поле\  $\mathbf{u}\in \mathbf{A}^{2k}$.
Рассмотрим главную часть скалярного произведения в $\mathbf{A}^{2k}(G)$:	
$$(\mathbf{u},\mathbf{w}_{\eta})_{2k}\equiv 
((\nabla \mathrm{div})^{k}\,\mathbf{u},(\nabla \mathrm{div})^{k}\,\mathbf{w}_{\eta}).$$
Проинтегрируем по частям:	
\[	(\mathbf{u},\mathbf{w}_{\eta})_{2k}=
((\nabla \mathrm{div})^{2k}\,\mathbf{u},\,\mathbf{w}_{\eta})=
\int_G
\mathbf{v}\cdot (\mathbf{w}_{\eta})\, d\mathbf{x}.   \eqno{(2.7)}\]
Левая часть имеет
предел при $\eta\rightarrow 0$, равный $(\mathbf{u},\mathbf{w})_{2k}$.
Следовательно,  правая часть также будет иметь предел и интеграл
$\int_G \mathbf{v}\cdot \mathbf{w}\, d\mathbf{x}$ существует при
любой $\mathbf{w}\in\mathbf{A}^{2k}(G)$. Кроме того,  из     
неравенства Коши-Буняковского следует оценка этого интеграла:
$$\left|\int_G
\mathbf{v}\cdot \mathbf{w}\, d\mathbf{x}\right|\leq
\|\mathbf{u}\|_{\mathbf{A}^{2k}}\|\mathbf{w}\|_{\mathbf{A}^{2k}}.   $$ 
Значит, $\mathbf{v}$ есть линейный функционал из
$(\mathbf{A}_{0}^{2k})^*$.
 
Применим его к полям $\mathbf{g}_i$,
составляющим базис пространства $\mathcal{A}_H(G)$.  Учитывая,
что  $\nabla \mathrm{div}\,\mathbf{g}_i=0$,
получим \[\int_G \mathbf{v}\cdot \mathbf{g}_i\, d\mathbf{x}=0. \quad
i=1,...,\rho_1.                  \eqno{(2.8)}  \] Итак,  на  полях
$\mathbf{w}$, отличающихся на вектор-функцию $\mathbf{g}$ из
$\mathcal{A}_H(G)$, его значения совпадают.
Пусть 	$\mathcal{A}/\mathcal{A}_H$ - фактор-пространство   	$\mathcal{A}(G)$ по $\mathcal{A}_H$   (пространство классов смежности). $\mathcal{A}^{2k}(G)=\mathbf{A}^{2k}(G)/\mathcal{A}_H$, его элементы имеют вид: $\mathbf{w}+\mathbf{g}$, где   $\nabla \mathrm{div}\,\mathbf{g}=0$.

\subsection{Оператор    $\mathcal{N}^{2k}_d$ в фактор-пространстве	 $\mathcal{A}^{2k}$ }
	\qquad    Шаг 2-й:
	
{\it	 Пространство
$\mathcal{A}^{2k}(G)$ становится гильбертовым, если  ввести скалярное произведение}
$$\{\mathbf{u},\mathbf{w}\}_{2k}\equiv (\mathbf{u},\mathbf{w})_{2k}=
((\nabla \mathrm{div})^{k}\,\mathbf{u},(\nabla \mathrm{div})^{k}\,\mathbf{w}).   \eqno{(2.9)}$$

Воспользуемся ортонормированным базисом в $\mathbf{A}^{0}$.
При $ \mathbf{f}\in \mathcal{A}^{2k}$, 
$\mathbf{g}_\eta\in	\mathcal{A}^{2k}_0$, 	в	терминах рядов Фурье  оно имеет вид:
\[\{\mathbf{f},\mathbf{g}_\eta\}_{2k}\equiv 
( \mathcal{N}_d^{k}\,\mathbf{f}, \mathcal{N}_d^{k}\,\mathbf{g}_\eta)=
\sum_{j=1}^{\infty}{\nu}_{j}^{4k}
[(\mathbf{f},\mathbf{q}_{j})(\mathbf{g}_\eta,\mathbf{q}_{j})
],  \eqno{(2.10)}\]  так как  
\[ \mathcal{N}_d^k\mathbf{f}=
\lim_{n\rightarrow\infty}
(\nabla\mathrm{div})^k\,(\mathbf{f}^n_{\mathcal{A}})= (-1)^k
\sum_{j=1}^{\infty}\nu^{2k}_j
(\mathbf{f},\mathbf{q}^{}_{j})\mathbf{q}_{j}(\mathbf{x}),   \eqno{(2.11)}\]

Для того чтобы функционал $\rho$ служил элементом
$( \mathcal{A}^{2k}_0)^*$, нужно, чтобы скалярное произведение
$({\rho}(\mathbf{x}),\mathbf{w}(\mathbf{x}))$ существовало при всех
$\mathbf{w}(\mathbf{x})\in  \mathcal{A}^{2k}$ и удовлетворяло
неравенству: $({\rho}(\mathbf{x}),\mathbf{w}(\mathbf{x}))\leq
M_{2k}\,\| \mathcal{N}_d^k\,\mathbf{w}\|$. 	
Мы имеем
$$({\rho},\mathbf{w})=\sum_{j=1}^{\infty}
[(\mathbf{\rho},\mathbf{q}_{j})(\mathbf{w},\mathbf{q}_{j})
] =  \sum_{j=1}^{\infty}
[(\mathbf{\rho},\mathbf{q}_{j}) /  \nu^{2k}_j ] [ \nu^{2k}_j (\mathbf{w},\mathbf{q}_{j})
] \leq 
M_{2k}\,\| \mathcal{N}_d^k\,\mathbf{w}\|,    \eqno{(2.12)} $$ где
\begin{equation*} \label
{karo 11}
M_{2k}= \{\sum_{j=1}^{\infty}{\nu}_{j}^{-4k}(\mathbf{\rho},
\mathbf{q}_{j})^2\}^{1/2}.  
\end{equation*}

Знак равенства при заданных $(\mathbf{\rho},\mathbf{q}_{j})$
достижим. 	Значит, имеет место
\begin{lemma}
   Условие
\[M_{2k}^2= \sum_{j=1}^{\infty}{\nu}_{j}^{-4k}(\mathbf{\rho},
\mathbf{q}_{j})^2  <\infty.  \eqno(2.13)\]
		необходимо и достаточно для принадлежности
	${\rho}(\mathbf{x})$ к $( \mathcal{A}^{2k}_0)^*$.
\end{lemma}

 Величина $M_{2k}$  есть норма функционала $\rho $ в
$( \mathcal{A}^{2k}_0)^*$, которая совпадает с нормой элемента
\[ \mathcal{N}_d^{-k}\,\mathbf{f}=	(-1)^k
\sum_{j=1}^{\infty}\nu^{-2k}_j
(\mathbf{f},\mathbf{q}_{j})\mathbf{q}_{j}(\mathbf{x}) \quad \text{при} \quad  \mathbf{f} \in  \mathbf{A}^{-2k}. \eqno{(2.14)}\]

Шаг 3-й: {\it   пространство $ (\mathbf{A}^{2k}_0)^*$  отождествим с пространством $ \mathbf{A}^{-2k}$,} 

Скалярное произведение в нем	определим  как 
\[\{\mathbf{u},\mathbf{w}\}_{-2k}=
( \mathcal{N}_d^{-k}\,\mathbf{u}, \mathcal{N}_d^{-k}\,\mathbf{w}), \eqno{(2.15)}\]      а
Лемму 2.1 переформулируем 
так

\begin{theorem}
		При заданном $\mathbf{v}\in \mathcal{A}^*$ и $k\geq 1$
	уравнение $(\nabla \mathrm{div})^{2k}\,\mathbf{u}=\mathbf{v}$  разрешимо в
	пространстве $\mathbf{A}^{2k}$ тогда и только тогда, когда
	$\mathbf{v}\in\mathbf{A}^{-2k}$.\quad
	Его решение   $\mathbf{u}= \mathcal{N}_d^{-2k}\mathbf{v}$
	в  фактор-пространстве	$\mathcal{A}/\mathcal{A}_H$   определяется однозначно.
\end{theorem}

Действительно, если   функционал 
$\mathbf{v}\in ( \mathbf{A}^{2k}_0)^*$, то его норма  $M_{2k}<\infty$
и он принадлежит $\mathbf{A}^{-2k}$, так как 	$\{\mathbf{v},\mathbf{v}\}_{-2k}=
( \mathcal{N}_d^{-k}\,\mathbf{v},  \mathcal{N}_d^{-k}\,\mathbf{v})=M_{2k}^2$.

Ряд $ \mathcal{N}_d^{k}\,\mathbf{u}= \mathcal{N}_d^{k}\,[  \mathcal{N}_d^{-2k},\mathbf{v}]=  \mathcal{N}_d^{-k}\,\mathbf{v}$ сходится в $\mathcal{A}_{\gamma}$, так как 
$( \mathcal{N}_d^{-k}\,\mathbf{v},  \mathcal{N}_d^{-k}\,\mathbf{v})=M_{2k}^2$.
Элемент $\mathbf{u}$ принадлежит  $\mathbf{A}^{2k}$ и удовлеторяет уравнению \newline  $(\nabla \mathrm{div})^{2k}\,\mathbf{u}= \mathcal{N}_d^{2k}\,[ \mathcal{N}_d^{-2k}\,\mathbf{v}]= \mathbf{v}$,   так как квадрат его нормы
\[\{\mathbf{u},\mathbf{u}\}_{2k}=
( \mathcal{N}_d^{k}\,\mathbf{u},  \mathcal{N}_d^{k}\,\mathbf{u})=( \mathcal{N}_d^{-k}\,\mathbf{v},  \mathcal{N}_d^{-k}\,\mathbf{v})=\{\mathbf{v},\mathbf{v}\}_{-m}=M_{2k}^2<\infty. \]
Однозначноcть решения вытекает из определения и обратимости операторов   $ \mathcal{N}_d$. Теорема доказана.

\section{Пространство  $\mathbf{V}^{0}$ вихревых полей }
  Joshida и Giga  \cite{yogi}  рассмотрели  пространство    $\mathbf{W}^1=\{\textbf{ f}\in  \mathbf{V}^0,\, \mathrm{rot}\,\mathbf{ f} \in\mathbf{V}^{0} \}$,\footnote { Пространства ${L}_{\sigma}^2(\Omega)$,   ${L}_{H}^2(\Omega)$, ${L}_{\Sigma}^2$,  ${H}_{\Sigma \,\Sigma}^1$  в  \cite{yogi} мы обозначаем  как ${\mathcal{{B}}}(\Omega)$,   $\mathcal{B}_H(\Omega)$,  $\mathbf{V}^{0}$ и $\mathbf{W}^1$.}
     и  оператор $S$ c    областью определения  $ \mathbf{W}^{1}$,   который совпадает с $\mathrm{rot}\mathbf{u}$, если   $\mathbf{u}\in  \mathbf{W}^{1}$.
 Они доказали, что   оператор $S$ самосопряжён в $\mathbf{V}^{0}$,
	спектр	 $\sigma(S)$  точечный и действительный  $\sigma(S)=\sigma_p(S)\subset\mathbb{R}$. 	
	
	Оператор  $S$ имеет компактный обратный из  $\mathbf{V}^{0}$ в   $ \mathbf{W}^{1} \subset \mathbf{H}^{1}$.
	 Значит,  пространство  $\mathbf{V}^{0}$,
	 составленное из  собственных полей оператора  $S$ имеет полный ортогональный базис.  Сам базис  авторы  \cite{yogi}       не рассматривали. 
	 
	 В  $\S 1$  мы рассмотрели  базис   $\mathbf{V}^{0}$,  его свойства в   области $G$,  поля $\mathbf{q}^{\pm}_{j}\in \mathbf{C}^\infty(\bar{G})$,   и	 их  формулы  (1.26)-(1.30))  в шаре  $B$   \cite{saUMJ13,saVS20}.

 Зафиксируем ортонормированный базис   $\{\mathbf{q}^{+}_{j}, \mathbf{q}^{-}_{j}\}$ в  $\mathbf{V}^{0}$ 

 \[\mathrm{rot}
\mathbf{q}_{j}^{\pm}=\pm\lambda_j\, \mathbf{q}_{j}^{\pm}, \quad    \gamma\mathbf{n}\cdot\mathbf{q}_{j}^{\pm}=0,   \quad    \|\mathbf{q}^{\pm}_{j}\|=1,   \quad  j=1, 2,  ... .      \eqno{(3.1)}  \]  
	
В этом базисе элементы    $\mathbf{V}^{0}(G)$ представляются рядами Фурье:		
\[\mathcal{P}_{\mathbf{V}}\,\mathbf{f}\equiv\mathbf{f}_{\mathbf{V}}=		
\lim_{n\rightarrow\infty}(\mathbf{f}^{n}_{\mathbf{V}})=
\sum_{j=1}^{\infty}
[(\mathbf{f},\mathbf{q}^{+}_{j})\mathbf{q}^{+}_{j}+
(\mathbf{f},\mathbf{q}^{-}_{j})\mathbf{q}^{-}_{j}],   \eqno{(3.2)}\]			
где $\mathbf{f}^{n}_{\mathbf{V}}$  - частичные суммы ряда,    а операторы   $S$  и  $S^{-1}$ -преобразованиями этих рядов: \[S\mathbf{f}_{\mathbf{V}}=\lim_{n\rightarrow\infty}\mathrm{rot}
(\mathbf{f}^n_{\mathbf{V}})=\sum_{j=1}^{\infty}\lambda_j
[(\mathbf{f},\mathbf{q}^{+}_{j})\mathbf{q}^{+}_{j}-
(\mathbf{f},\mathbf{q}^{-}_{j})\mathbf{q}^{-}_{j}],    \eqno{(3.3)}\]		
\begin{equation*}\klabel{obr__2_}S^{-1}\mathbf{f}_{\mathbf{V}}=
\sum_{j=1}^{\infty}\lambda_j^{-1}
[(\mathbf{f},\mathbf{q}^{+}_{j})\mathbf{q}^{+}_{j}-
(\mathbf{f},\mathbf{q}^{-}_{j})\mathbf{q}^{-}_{j}]. \eqno{(3.4)} \end{equation*}
\subsection{Подпространста $\mathbf{W}^m$ в $\mathbf{V}^0$}	В \cite{saVS20} мы ввели пространства	{\footnote{ \small    Авторы   \cite{yogi} рассматривали только пространство $H^1_{\Sigma \Sigma}$.}}							
$$\mathbf{W}^m=\{\textbf{ f}\in  \mathbf{V}^0,...,( \mathrm{rot})^m\mathbf{ f}\in \mathbf{V}^{0} \},   \quad m =1, 2, 3, ... \eqno{(3.5)}$$	
Согласно  оценке (1.12)  пространство	
$\mathbf{W}^m\subset \mathbf{H}^m$  и 
$\mathbf{W}^m$   еcть  проекция $ \mathbf{H}^m$  на  $\mathbf{V}^0$, так как   для любого поля   $\mathbf{f}\in  \mathbf{H}^m $  его проекция $\mathcal{P}_{V}\mathbf{f}\in  \mathbf{W}^{m}$; если  же $\mathbf{f}_{\mathbf V}\in   \mathbf{W}^{m}$,   то $\mathcal{P}_{V}\mathbf{f}_{\mathbf{V}}= \mathbf{f}_{\mathbf{V}}$,  а 	его проекция на 	  	  $\mathcal{A}$ равна 0.  

Замыкание в $\mathbf{W}^{m}(G)$ пространства
$\mathbf{C}^{\infty}_0(G)$ обозначается через
$\mathbf{W}^{m}_0$,  а пространство    $(\mathbf{W}^{m}_0)^*$ сопряжено с $\mathbf{W}_0^{m}$, - это пространство линейных непрерывных функционалов над  $\mathbf{W}^{m}_0$.  Они  равны  нулю на пространстве  $\mathcal{B}_{H}$. 
Объединение  пространств   $(\mathbf{W}^{k}_0)^*$  обозначим $\mathcal{B}^*$.	

Пусть $m\geq 1$, а область $G$ такова, что $\rho= dim \mathcal{B}_H>0$.	 Оператор     $S\mathbf{u}$ совпадает с $\mathrm{rot}\mathbf{u}$, если   $\mathbf{u}\in  \mathbf{W}^{1}\equiv \mathcal{D}(S) $. Поэтому оператор $(\mathrm{rot})^{2m}$ на $\mathbf{W}^m\subset \mathbf{W}^1$ совпадает с $S^{2m}$. 
\subsection{Оператор  $S^{2m}$  в пространстве $\mathbf{W}^m$} Основное утверждение:

{\it Оператор  $S^{2m}$ отображает  пространство		$\mathbf{W}^m$ на $\mathbf{W}^{-m}$ и обратно.}

Этапы доказательства:

Шаг 1-й:  {\it  Оператор $S^{2m}$ отображает  пространство
	$\mathbf{W}^m(G)$ на $(\mathbf{W}^{m}_0)^*$ .} 

		Действительно,  	пусть $\mathbf{w}$ произвольный элемент из $\mathbf{W}^m(G)$, а
$\mathbf{w}_{\eta}$ -- средняя вектор-функция для него,
$\mathbf{w}_{\eta}\in\mathbf{W}^m_0$, поле  $\mathbf{u}\in\mathbf{W}^m$.
 Рассмотрим главную часть скалярного произведения в $\mathbf{W}^m$:	
$(\mathbf{u},\mathbf{w}_{\eta})_m\equiv 
(\mathrm{rot}^m\,\mathbf{u},\mathrm{rot}^m\,\mathbf{w}_{\eta})$.

Проинтегрируем по частям:	
\[	(\mathbf{u},\mathbf{w}_{\eta})_m\equiv	(\mathrm{rot}^{2m}\,\mathbf{u},\,\mathbf{w}_{\eta})=
\int_G
\mathbf{v}\cdot (\mathbf{w}_{\eta})\, d\mathbf{x}.  \eqno{(3.6)}\]
Левая часть имеет
предел при $\eta\rightarrow 0$, равный $(\mathbf{u},\mathbf{w})_m$.
Следовательно,  правая часть также будет иметь предел и интеграл
$\int_G \mathbf{v}\cdot \mathbf{w}\, d\mathbf{x}$ существует при
любой $\mathbf{w}\in\mathbf{W}^m(G)$. Кроме того,  из     
неравенства Коши-Буняковского следует оценка этого интеграла:
$$\left|\int_G
\mathbf{v}\cdot \mathbf{w}\, d\mathbf{x}\right|\leq
\|\mathbf{u}\|_{\mathbf{W}^m}\|\mathbf{w}\|_{\mathbf{W}^m}.    \eqno{(3.7)}$$ 
Значит, $\mathbf{v}$ есть линейный функционал из
$(\mathbf{W}_0^m)^*$.

Применим его к полям $\mathbf{h}_i$,
составляющим базис пространства $\mathcal{B}_H(G)$.  Учитывая,
что  $\mathrm{rot}\,\mathbf{h}_i=0$,
получим \[\int_G \mathbf{v}\cdot \mathbf{h}_i\, d\mathbf{x}=0. \quad
i=1,...,\rho.    \eqno{(3.8)}\] Итак,  на  полях
$\mathbf{w}$, отличающихся на вектор-функцию $\mathbf{h}$ из
$\mathcal{B}_H(G)$, его значения совпадают.

Пусть 	$\mathcal{B}/\mathcal{B}_H$ - фактор-пространство   	$\mathcal{B}(G)$ по $\mathcal{B}_H$   (пространство классов смежности)
$\mathcal{B}^m(G)=\mathbf{W}^m(G)/\mathcal{B}_H$, его элементы имеют вид: $\mathbf{w}+\mathbf{h}$,   где  $\mathrm{rot}\, \mathbf{h} = 0$.
	\subsection { Оператор $S^{2m}$ в фактор-пространстве $\mathcal{B}^m$}   Шаг 2-й.
	
	 {\it  Пространство
$\mathcal{B}^m(G)$ становится гильбертовым, если на
нем ввести скалярное произведение $\{\mathbf{u},\mathbf{w}\}_m\equiv 	(\mathrm{rot}^m\,\mathbf{u},\mathrm{rot}^m\,\mathbf{w})$.	}

Воспользуемся ортонормированным базисом в $\mathbf{V}^0$.

При $ \mathbf{f}\in \mathcal{B}^m$, 
$\mathbf{g}_\eta\in	\mathcal{B}^m_0$, 	в	терминах рядов Фурье  оно выражается так:
\[\{\mathbf{f},\mathbf{g}_\eta\}_m\equiv 
(S^m\,\mathbf{f},S^m\,\mathbf{g}_\eta)=
\sum_{j=1}^{\infty}({\lambda}_{j}^{2m})
[(\mathbf{f},\mathbf{q}_{j}^+)(\mathbf{g}_\eta,\mathbf{q}_{j}^+)
+(\mathbf{f},\mathbf{q}_{j}^-)(\mathbf{g}_\eta,\mathbf{q}_{j}^-)],  \eqno{(3.9)}\]  так как  \[S^m\mathbf{f}=  
\lim_{n\rightarrow\infty}\mathrm{rot}^m
(\mathbf{f}^n_{\mathbf{V}})=	
 \sum_{j=1}^{\infty}\lambda_j^m
[(\mathbf{f},\mathbf{q}^{+}_{j})\mathbf{q}^{+}_{j}+ (-1)^m
(\mathbf{f},\mathbf{q}^{-}_{j})\mathbf{q}^{-}_{j}].   \eqno{(3.10)}\]

Для того чтобы функционал $\rho$ служил элементом
$(\mathcal{B}^m_0)^*$, нужно, чтобы скалярное произведение
$({\rho}(\mathbf{x}),\mathbf{w}(\mathbf{x}))$ существовало при всех
$\mathbf{w}(\mathbf{x})\in \mathcal{B}^m$ и удовлетворяло
неравенству: $({\rho}(\mathbf{x}),\mathbf{w}(\mathbf{x}))\leq
K_m\,\|S^m\,\mathbf{w}\|$. 	
Мы имеем
$$({\rho},\mathbf{w})=\sum_{j=1}^{\infty}
[(\mathbf{\rho},\mathbf{q}_{j}^+)(\mathbf{w},\mathbf{q}_{j}^+)
+(\mathbf{\rho},\mathbf{q}_{j}^-)(\mathbf{w},\mathbf{q}_{j}^-) ]\leq 
K_m\,\|S^m\,\mathbf{w}\|,    \eqno{(3.11)} $$ где
\begin{equation*} \label
{karo 1}
K_m= \{\sum_{j=1}^{\infty}({\lambda}_{j}^{-2m})[(\mathbf{\rho},
\mathbf{q}^{+}_{j})^2+(\mathbf{\rho},\mathbf{q}^{-}_{j})^2]\}^{1/2}.  
\end{equation*}
%\times \{\sum_{j=1}^{\infty}({\lambda}_{j}^{2})[(\mathbf{w},
%\mathbf{q}^{+}_{j})^2+(\mathbf{w},\mathbf{q}^{-}_{j})^2]\}^{1/2}.\]
Знак равенства при заданных $(\mathbf{\rho},\mathbf{q}^{\pm}_{j})$
достижим. 	Значит, имеет место

	\begin{lemma}Условие
\begin{equation*}
K_m^2=\sum_{j=1}^{\infty}({\lambda}_{j}^{-2m})[(\mathbf{\rho},
\mathbf{q}^{+}_{j})^2+(\mathbf{\rho},\mathbf{q}^{-}_{j})^2]<\infty.  \eqno{(3.12)}
\end{equation*} %\label{spb 1}
необходимо и достаточно для принадлежности
${\rho}(\mathbf{x})$ к $(\mathcal{B}_0^m)^*$.
\end{lemma}

Величина $K_m$  есть норма функционала $\rho $ в
$(\mathcal{B}^m_0)^*$, которая совпадает с нормой элемента
\[S^{-m}\,\mathbf{f}=	
\sum_{j=1}^{\infty}\lambda_j^{-m}
[(\mathbf{f},\mathbf{q}^{+}_{j})\mathbf{q}^{+}_{j}+ (-1)^m
(\mathbf{f},\mathbf{q}^{-}_{j})\mathbf{q}^{-}_{j}]
 \quad \text{при} \quad  \mathbf{f} \in  \mathbf{W}^{-m}.     \eqno{(3.13)}\]

%$\|S^{-m}\,\mathbf{\rho}\|=K_m$, при 	$\rho\in \mathbf{W}^{-m}$. 

\subsection {Пространства $(\mathcal{B}^m_0)^*$  и $ \mathbf{W}^{-m}$}  Шаг 3-й.

{\it   Пространство $(\mathcal{B}_0^m)^*$  отождествим с
пространством $\mathbf{W}^{-m}$.}

Скалярное произведение в нем	определим как 
\[\{\mathbf{u},\mathbf{w}\}_{-m}=
(S^{-m}\,\mathbf{u},S^{-m}\,\mathbf{w}), \eqno{(3.14)}\]      а
Лемму 3.1 переформулируем так.

\begin{theorem} 
	При заданном $\mathbf{v}\in \mathcal{B}^*$ и $m\geq 1$
	уравнение $\mathrm{rot}^{2m}\,\mathbf{u}=\mathbf{v}$  разрешимо в
	пространстве $\mathcal{B}^{m}(G)$ тогда и только тогда, когда
	$\mathbf{v}\in\mathbf{W}^{-m}(G)$.\quad
	Его решение   $\mathbf{u}=S^{-2m}\mathbf{v}$ 	в классе   $\mathcal{B}(G)/\mathcal{B}_H$  определяется однозначно.
	\end{theorem}

	 Действительно, если   функционал 
$\mathbf{v}\in (\mathcal{B}_0^m)^*$, то его норма  $K_m<\infty$
и он принадлежит $\mathbf{W}^{-m}$, так как 	$\{\mathbf{v},\mathbf{v}\}_{-m}=
(S^{-m}\,\mathbf{v}, S^{-m}\,\mathbf{v})=K_m^2$.

Ряд $S^{m}\,\mathbf{u}=S^{m}\,[ S^{-2m}\,\mathbf{v}]= S^{-m}\,\mathbf{v}$ сходится в $\mathbf{V}^0$, так как 
$(S^{-m}\,\mathbf{v}, S^{-m}\,\mathbf{v})=K_m^2$.	
Элемент $\mathbf{u}$ принадлежит  $\mathcal{B}^m$ и удовлеторяет уравнению $\mathrm{rot}^{2m}\,\mathbf{u}=S^{2m}\,[ S^{-2m}\,\mathbf{v}]= \mathbf{v}$ ,  так как квадрат его нормы
$\{\mathbf{u},\mathbf{u}\}_{m}=
(S^{m}\,\mathbf{u}, S^{m}\,\mathbf{u})=(S^{-m}\,\mathbf{v}, S^{-m}\,\mathbf{v})=\{\mathbf{v},\mathbf{v}\}_{-m}=K_m^2<\infty$. 
Однозначноcть решения вытекает из определения и обратимости операторов   $S$. Теорема доказана.

Цепь  вложений  пространств $ \mathbf{W}^{k}$ имеет вид:
\[\subset \mathbf{W}^{k}\subset...\subset \mathbf{W}^1\subset \mathbf{V}^0\subset  \mathbf{W}^{-1}\subset...\subset \mathbf{W}^{-k}\subset. \eqno{(3.15)} \]	

Теорема  4 показывает, что в этой цепи имеется симметрия:
{\it для любого $m\geq 1$  отображения 
	$S^{2m}: \mathbf{W}^m\to \mathbf{W}^{-m}$ и
$S^{-2m}:\mathbf{W}^{-m}\to \mathbf{W}^m$ взаимно обратны  и  	} $$\{\mathbf{u},\mathbf{u}\}_{m}=
	(S^{m}\,\mathbf{u}, S^{m}\,\mathbf{u})=(S^{-m}\,\mathbf{v}, S^{-m}\,\mathbf{v})=\{\mathbf{v},\mathbf{v}\}_{-m}.    \eqno{(3.16)}$$ 

Отображения $S:\mathbf{W}^{m+1} \rightarrow \mathbf{W}^{m}$ при $m\geq 1$ и  	$S^{-1}:\mathbf{W}^{m} \to \mathbf{W}^{m+1}$  также взаимно обратны.  
Причем 
$ \|S^{-1}\mathbf{f}\|^2_{\mathbf{W}^{m+1}} 
\leq  C^2_m \|\mathbf{f}\|^2_{\mathbf{W}^{m}},$ где % \newline
$C^2_m= max _j(1+1/{\lambda}_{j}^{2m})$   и $1/{\lambda}_{j}\to 0$ при $j\to \infty $, а
$\|  S\mathbf{f}\|^2_{ \mathbf{W}^m}\leq  
C^{-2}_m \|\mathbf{f}\|^2_{\mathbf{W}^{m+1}}$. 

\subsection {Операторы $(S+\lambda\,I)$ в пространствах $ \mathbf{W}^{m}$}
Автор  рассмотрел также  операторы  $(S+\lambda\,I):\mathbf{W}^{m+1} \rightarrow \mathbf{W}^{m}$ при $m\geq 0$  и   $(S+\lambda\,I)^{-1}:\mathbf{W}^{m} \to \mathbf{W}^{m+1}$  и  доказал  в \cite{saVS20} их разрешимость по Фредгольму:
\begin{theorem}\klabel{S  1}
	Оператор $S+\lambda I:
	\mathbf{W}^{m+1}(G)\rightarrow \mathbf{W}^{m}(G)$ непрерывен и
	однозначно обратим, если $\lambda$ не принадлежит спектру
	оператора $S$. Его обратный задается формулой	
	\begin{equation*}\klabel{sp__3_}
	(S+\lambda I)^{-1}\mathbf{f}=
	\sum_{j=1}^{\infty}\left[
	\frac{(\mathbf{f},{\mathbf{q}}_{j}^+)}{\lambda+\lambda_{j}}
	\mathbf{q}_{j}^+(\mathbf{x})+
	\frac{(\mathbf{f},{\mathbf{q}}_{j}^-)}{\lambda-\lambda_{j}}
	\mathbf{q}_{j}^-(\mathbf{x})\right]. \eqno (3.17) \end{equation*}
	Если же $\lambda=\lambda_{j_0}$, то он обратим тогда и только
	тогда, когда
	\begin{equation*}
	\klabel{urz 1_}\int_G \mathbf{f}\cdot \mathbf{q_j^-} d\mathbf{x}=0\quad
	\text{для}\quad\forall \mathbf{q_j^-}: \lambda_j=\lambda_{j_0}. \eqno{(3.18)}
	\end{equation*}
	Ядро оператора $S+\lambda_{j_0} I$ конечномерно и определяется
	собственными функциями $\mathbf{q_j^-}(\mathbf{x})$, собственные
	значения которых равны $\lambda_{j_0}$:
	\begin{equation*} \klabel{ker__1_}
	Ker(S+\lambda_{j_0} I)= \sum_{\lambda_j=\lambda_{j_0}}
	c_j \mathbf{q}^{-}_{j}(\mathbf{x})\quad\forall \,c_j\in
	\mathbb{R}.\eqno{(3.19)} \end{equation*}\end{theorem}
В \cite{saVS20} доказаны  также  оценки  
\[ \|  (S+\lambda I)\mathbf{f}\|^2_{ \mathbf{W}^m}\leq  c^2_m \|\mathbf{f}\|^2_{\mathbf{W}^{m+1}}, \,\,  \|(S+\lambda I)^{-1}\mathbf{f}\|^2_{\mathbf{W}^{m+1}}\leq  C^2_m \|\mathbf{f}\|^2_{\mathbf{W}^{m}},  \eqno{(3.20 )} \]
где  постоянные $  c^2_m = max _j \, (a_{m,j}^+, a_{m,j}^-)$ \,  и $ a_{m,j}^{\pm }= (|1\pm \lambda/\lambda_j|^2/(1+1/{\lambda}_{j}^{2m})) $; 
\[C^2_m= max _j \, (A_{m,j}^+, A_{m,j}j^-) \,\,\text{ и } A_{m,j}^{\pm }= ((1+1/{\lambda}_{j}^{2m}) /|1\pm \lambda/\lambda_j|^2.  \eqno{(3.21 )} \] 
числа $ a_{m,j}^{\pm }$  и $ A_{m,j}^{\pm }<\infty$ (при больших $\lambda_{j}$ они находятся в окрестности 1).

%end{documen}
%\newpage 

\subsection {Соотношения между пространствами  $ \mathbf{W}^{k}(\Omega)$,    $ \mathbf{H}^{k}(\Omega)$	и    $ \mathbf{C}^{k-2}(\bar{ \Omega})$}

В области $\Omega$,  гомеоморфной шару,  в \cite{saVS20}  доказана
\begin{theorem}\klabel{W}
	Для того, чтобы $\mathbf{f}\in \mathbf{V}^0(\Omega)$ разлагалась в ряд Фурье
	\begin{equation*} \klabel{rof 1}
	\mathbf{f}(\mathbf{x})=\sum_{j=1}^{\infty}
	((\mathbf{f},\mathbf{q}_{j}^+)\mathbf{q}_{j}^+(\mathbf{x})
	+(\mathbf{f},\mathbf{q}_{j}^-)\mathbf{q}_{j}^-(\mathbf{x})), \quad
	\|\mathbf{q}_{j}^{\pm}\| =1,\eqno{(3.22)}
	\end{equation*}
	по собственным вектор-функциям $\mathbf{q}_{j}^{\pm}(\mathbf{x})$
	ротора в области $\Omega$,
	сходящийся в норме
	пространства Соболева $\mathbf{H}^k(\Omega)$, необходимо и достаточно,
	чтобы $\mathbf{f}$ принадлежала $\mathbf{W}^k(\Omega)$.
	
	Если $\mathbf{f}\in \mathbf{W}^k(\Omega)$,
	то существует такая постоянная $C>0$, не зависящая от
	$\mathbf{f}$, что
	\begin{equation*} \klabel{orf 3}
	\sum_{j} {\lambda}_{j}^{2k} ((\mathbf{f},\mathbf{q}_{j}^+)^2
	+(\mathbf{f},\mathbf{q}_{j}^-)^2)\leq C
	\|\mathbf{f}\|^2_{\mathbf{H}^k(\Omega)}.\eqno{(3.23)}
	\end{equation*}
	Если $k\geq 2$, то  вектор-функция $\mathbf{f}$ из	$\mathbf{W}^k(\Omega)$
	разлагается в в ряд (3.22), сходящийся в пространстве 
	$\mathbf{C}^{k-2}(\overline{\Omega})$.     \end{theorem}

	\begin{corol} Вектор-функция $\mathbf{f}\in
	\mathbf{V}^0\cap\mathbf{C}^{\infty}_0(\Omega)$
	разлагается в	ряд	(3.22),
	сходящийся в любом из пространств
	$\mathbf{C}^{k}(\overline{\Omega})$, \, $ k\in \mathbb{N}$.\end{corol}

Эти результаты дополняют известные в теории рядов Фурье утверждения
(см. Теорема 7 в $\S 4$ гл. 2 в [5], Теорема 8 в $\S 2$ гл. 4 в [2]).

Таким образом, $\mathbf{W}^k(\Omega)$ играют роль
пространств Соболева в классе вихревых  полей. 

\section{Семейство пространств Соболева: классы  $ \mathbf{C}(2k,m)\equiv \mathbf{A}^{2k} \oplus \mathbf{W}^m$}
	Пусть   область $\Omega$   гомеоморфна шару,  а    $ \mathbf{C}(2k,m)$-  классы в cети пространств Соболева $\mathbf{A}^{2k} \oplus \mathbf{W}^m$,   числа $k, m $  -целые.
\subsection{Операторы $ \mathcal{N}^p_d$ и $S^p$ в шкалах пространств   $ \mathbf{A}^{2k}(\Omega)$  и  $ \mathbf{W}^m(\Omega)$} 
 Если собственные поля $ \mathbf{q}_{j}(\mathbf{x})$  и  $\mathbf{q}_{j}^\pm(\mathbf{x})$  градиента дивергенции и ротора известны, то  элементы    $\mathbf{f}_\mathcal{A}\in \mathcal{A}$   и  $\mathbf{f}_\mathbf{V}\in \mathcal{B}=\mathbf{V}^0 $   
представляются   рядами Фурье: 
\[\mathbf{f}_\mathcal{A}=
\sum_{j=1}^{\infty} (\mathbf{f},{\mathbf{q}}_{j}) \mathbf{q}_{j }(\mathbf{x}), \quad 	\mathbf{f}_\mathbf{V}=\sum_{j=1}^{\infty}\left[(\mathbf{f},{\mathbf{q}}_{j}^{+})
\mathbf{q}_{j}^{+}(\mathbf{x})+
(\mathbf{f},{\mathbf{q}}_{j}^{-})
\mathbf{q}_{j}^{-}(\mathbf{x})\right],     \eqno{(4.1)}\]
а   элементы    $\mathbf{f}$  из $\mathbf{L}_2(\Omega)$ -- их суммой $\mathbf{f}_\mathcal{A}+\mathbf{f}_\mathbf{V}$.  

 Причём $  \text{rot}\,\mathbf{f}_\mathcal{A}=0$ и $\text{div}\,\mathbf{f}_\mathbf{V}=0$, поэтому  $ \text{div}\,\mathbf{f}=  \text{div}\,\mathbf{f}_\mathcal{A}$,  а   $\text{rot}\,\mathbf{f}=  \text{rot}\,\mathbf{f}_\mathbf{V}$, 

Скалярное произведение  $(\mathbf{f}, \mathbf{g})$ полей
$\mathbf{f}$  и $ \mathbf{g}$  из $\mathbf{L}_2(\Omega)$ равно   $(\mathbf{f}_\mathcal{A}, \mathbf{g}_\mathcal{A})+(\mathbf{f}_\mathbf{V}, \mathbf{g}_\mathbf{V})$.

Действия операторов    $\mathcal{N}_d^p$  в $ \mathcal{A}$, \,   $S^p$ в  	
$ \mathcal{B}$  и  обратных при $p=1, 2, ...$  имеют вид
\[  \mathcal{N}^p_d \mathbf{f}_\mathcal{A}=(-1)^p\sum_{j=1}^{\infty}(\nu^{2p}_j) (\mathbf{f},{\mathbf{q}}_{j}) \mathbf{q}_{j }, \quad  	S^p \mathbf{f}_\mathbf{V}=\sum_{j=1}^{\infty}\lambda^p_{j}\left[(\mathbf{f},{\mathbf{q}}_{j}^{+})
\mathbf{q}_{j}^{+}+(-1)^p(\mathbf{f},{\mathbf{q}}_{j}^{-})
\mathbf{q}_{j}^{-}\right],   \]
\[  \mathcal{N}_d  ^{-p}\mathbf{f}_\mathcal{A}=(-1)^p
\sum_{j=1}^{\infty}(\nu^{-2p}_j) (\mathbf{f},{\mathbf{q}}_{j}) \mathbf{q}_{j }, \quad    	S^{-p}\mathbf{f}_\mathbf{V}=\sum_{j=1}^{\infty}\lambda_{j}^{-p}\left[(\mathbf{f},{\mathbf{q}}_{j}^{+}) \mathbf{q}_{j}^{+}+(-1)^p
(\mathbf{f},{\mathbf{q}}_{j}^{-}) \mathbf{q}_{j}^{-} \right].     \]

Так как $ \mathcal{A}_{H}= \mathcal{B}_{H}=  \emptyset$,  то  $ \mathcal{A}_{\gamma }^0= \mathcal{A}_{\gamma }(B)=\{\nabla\,h, h\in H^2(B): \gamma (\mathbf{n}\cdot \nabla) h =0 \}$,
$ \mathbf{V}^0=\{\mathbf{g}\in \mathbf{L}_2(B), \,\mathbf{g} \perp \mathcal{A}, \,\, \text{div}\,\mathbf{g}=0, \,\gamma (\mathbf{n}\cdot \mathbf{g})  =0 \}$,  и  пространства
\[ \mathbf{A}^{2k}\equiv \{\mathbf{f}\in \mathcal{A}_{\gamma }, ...,(\nabla  \mathrm{div})^k \mathbf{f}\in \mathcal{A}_{\gamma }\}\quad \text{и} \quad  \mathbf{W}^m \equiv \{\mathbf{g}\in \mathbf{V}^0,..., (\mathrm{rot})^m \mathbf{g} \in \mathbf{V}^0\} \eqno{(4.2)} \]  при  $k\geq 1, m\geq1$; \quad 
$ \mathbf{A}^{0}\equiv  \mathcal{A}_{\gamma }, \quad    \mathbf{W}^{0}\equiv  \mathbf{V}^{0}\equiv   \mathcal{B}$. 

Имеют место вложения
\[...\subset   \mathbf{A}^{2k}\subset...   \subset  \mathbf{A}^{2} \subset   \mathbf{A}^{0} \subset  \mathbf{A}^{-2} \subset ... \subset \mathbf{A}^{-2k}\subset ...,   \eqno{(4.3)} \]
%\quad 
\[...\subset \mathbf{W}^m\subset ...\subset \mathbf{W}^1\subset \mathbf{V}^{0}\subset \mathbf{W}^{-1}\subset...\subset \mathbf{W}^{-m}\subset.... \eqno{(4.4)}\]
В  шкале  пространств (4.3) согласно (14) $\S 1$  оператор   $\mathcal{N}_d$  действуют слева направо,  а  оператор   $\mathcal{N}^{-1}_d$  -cправа налево.

Операторы $S$ согласно (15) $\S 1$ действуют аналогично.

Прямые суммы пространств  $ \mathbf{A}^{2k}$  и 
$\mathbf{W}^m$  мы  обозначили  как $  \mathbf{C}(2k,m)$, они принадлежат   $\mathbf{L}_{2}(\Omega)$, если $k\geq 0, m\geq 0 $  -целые.

Операторы $(\mathcal{N}_d^{-1}, I)$ и $(I, S^{-1})$  и отображают  класс  $  \mathbf{C}(2k,m)$ на   $ \mathbf{C}(2(k+1),m)$ и   $ \mathbf{C}(2k,m+1)$,  а оператор  $(\mathcal{N}_d^{-p}, S^{-q})$ - на  класс  $ \mathbf{C}(2(k+p), m+q)$ при  $p, q>0$.

% $ \mathbf{C}$
Отметим,  что  операторы   $(\nabla \mathrm {div})^p$ и $(\mathrm {rot})^{2p}$, где  $ p$  - натуральное  число,
- аналоги  полигармонических  операторов $\Delta^p$
в классах $\mathcal{A}$   и  $\mathcal{B}$. 
\subsection {Краевые задачи в сети пространств 
		$ \mathbf{C}(2k, m)$} Классы  $ \mathbf{C}(2k, m)$ принадлежат двух-параметрическому семейству пространств Соболева.
	
 Пространство $\mathbf{A}^{2k}(\Omega)\subset \mathbf{H}^{2k}(\Omega)$ 
и является проекцией $\mathbf{H}^{2k}$ на $\mathcal{A}$,  а
$\mathbf{W}^{m}(\Omega)\subset \mathbf{H}^{m}(\Omega)$  
является проекцией $\mathbf{H}^{m}$ на $\mathcal{B}$, поэтому
 класс  $ \mathbf{C}(2k, 2k)$ совпадает с пространством Соболева    $\mathbf{H}^{2k}(\Omega)$.
%NSS (the Net of Sobolev Spaces). Рассмотрим задачи:
 \begin{problem}p.  Задано поле  $\mathbf{f}\in
	 \mathbf{C}(2k, m) \subset \mathbf{L}_{2}(\Omega) $.
	  Найти поле 
	$\mathbf{u}$  в $\mathbf{L}_{2}(\Omega)$ такое, что
	\begin{equation*}	\klabel{rot__4_}
{(\mathrm {rot})^{p}}\, \mathbf{u}+\lambda \mathbf{u}=
\mathbf{f} \quad \text{в}\quad \Omega,  \quad \gamma(\mathbf{n}\cdot \mathbf{u})=0.\eqno{(4.5)}\end{equation*}
 \end{problem}
{\it Это означает:} найти  поле $\mathbf{u}$  в $	\mathbf{L}_{2}(\Omega)$,  для которого  $(\mathbf{u},(\mathrm {rot}^{p}+\lambda\, I)\mathbf{ v})= (	\mathbf{f},	\mathbf{v})$ для любого поля $\mathbf{v}\in \,\mathbf{C}_0^{\infty}(\Omega) $, и 
 $ \gamma(\mathbf{n}\cdot \mathbf{u})=0$,  если след $\gamma(\mathbf{n}\cdot \mathbf{u})$ на границе $\Omega$  существует.	
 \begin{problem} \klabel{2} p. Задано поле  $\mathbf{f}\in
 	 \mathbf{C}(2k, m)  \subset \mathbf{L}_{2}(\Omega)$.  Найти поле 
	$\mathbf{w}$ в $\mathbf{L}_{2}(\Omega)$ такое, что
\begin{equation*}	\klabel{nd_ 4_}
(\nabla \mathrm {div})^p\, \mathbf{w}+\lambda \mathbf{w}=
\mathbf{f} \quad \text{в}\quad \Omega   \quad \text{и} \quad \gamma(\mathbf{n}\cdot \mathbf{w})=0,   \eqno{(4.6)}\end{equation*}
\end{problem}
если след $\gamma(\mathbf{n}\cdot \mathbf{w})$ на границе $\Omega$ существует.

Применим метод ортогонального проектирования  уравнений этих  задач в   пространстве $\mathbf{L}_{2}(\Omega)=  \mathcal{A} \oplus \mathcal{B} $ на подпространства $\mathcal{A}$ и $ \mathcal{B}$.
Рассмотрим   объемлющее пространство $\mathbf{L}_{2}(\Omega)$. 
Используя разложение полей   $\mathbf{f}$ ,   $\mathbf{u}$ и  $\mathbf{w}$  
из $\mathbf{L}_2(\Omega)$ -- в  суммы:  $\mathbf{f}_\mathcal{A}+\mathbf{f}_\mathbf{V}$,   $\mathbf{u}_\mathcal{A}+\mathbf{u}_\mathbf{V}$  и $\mathbf{w}_\mathcal{A}+\mathbf{w}_\mathbf{V}$ и 
расширения  $S$ и $\mathcal{N}_d$ операторов ротор и  градиент дивергенции, эти уравнения запишем  в виде  уравнений - проекций
 на $\mathcal{A}$ и  $\mathcal{B}$:
\begin{equation*}	\klabel{pro__4_}
\lambda\mathbf{u}_\mathcal{A}=\mathbf{f}_\mathcal{A},  \quad 
(S^{p}+\lambda I)\mathbf{u}_\mathbf{V}=\mathbf{f}_\mathbf{V},
\eqno{(4.7)}\end{equation*}
\begin{equation*}	\klabel{prond__4_}
(\mathcal{N}_d^{p}+\lambda  I)\mathbf{w}_\mathcal{A}=\mathbf{f}_\mathcal{A},
\quad
\lambda \mathbf{w}_\mathbf{V}=\mathbf{f}_\mathbf{V}, 
 \eqno{(4.8)}\end{equation*}
учитывая, что   $\text{rot} \,\mathbf{u}_\mathcal{A}=0$ в  $ \mathcal{A}$,   
$  \nabla\text{div}\,\mathbf{w}_\mathbf{V}=0$  в $ \mathcal{B}\equiv \mathbf{V}^0$.

 \begin{problem} \klabel{3} Задано поле  $\mathbf{f}\in
	 \mathbf{C}(2k, m) \subset \mathbf{L}_{2}(\Omega)$.  Найти поле 
	$\mathbf{u}$ в $\mathbf{L}_{2}(\Omega)$ такое, что
	\begin{equation*}	\klabel{nd_ 4_}
	\nabla \mathbf {div}\, \mathbf{u}+{\bf rot}\, \mathbf{u}+\lambda \mathbf{u}=
	\mathbf{f} \quad \text{в}\quad \Omega   \quad \text{и} \quad \gamma(\mathbf{n}\cdot \mathbf{u})=0,   \eqno{(4.9)}\end{equation*}
\end{problem}
если след $\gamma(\mathbf{n}\cdot \mathbf{u})$ на границе $\Omega$ существует.

Это уравнение эквивалентно двум   уравнениям-проекциям:
\[   (\mathcal{N}_d+\lambda I)\mathbf{u}_\mathcal{A}=\mathbf{f}_\mathcal{A}, 
\quad (S+\lambda I)\mathbf{u}_\mathbf{V}=\mathbf{f}_\mathbf{V}.  \eqno{(4.10)} \]

{\it Замечание.} Если пространство    $\mathcal{B}_H(G)$ не пусто и     $\lambda\neq 0$,  то уравнение  \newline  $(\nabla\mathrm{div} +\mathrm{rot} +\lambda  I)\mathbf{u}=\mathbf{f}$ распадается на три проекции
\[   (\mathcal{N}_d+\lambda I)\mathbf{u}_\mathcal{A}=\mathbf{f}_\mathcal{A}, 
\quad (S+\lambda I)\mathbf{u}_\mathbf{V}=\mathbf{f}_\mathbf{V}, \quad  
\lambda\mathbf{u}_{B_H}=\mathbf{f}_{B_H} \eqno{(4.11)} \]
-  уравнения  второго, первого и нулевого порядков, соответственно.

Согласно Теоремам 3 и 6 уравнения 
 $(\mathcal{N}_d+\lambda  I)\mathbf{u}_\mathcal{A}=\mathbf{f}_\mathcal{A}$  и 
 $ (S+\lambda I)\mathbf{u}_\mathbf{V}=\mathbf{f}_\mathbf{V}$  

разрешимы по Фредгольму.  В работе  \cite{saVS20}
  доказана

\begin{theorem}\klabel{1 RG}
 При $\lambda\neq Sp (rot)$ единственное решение задачи 1.1 имеет вид:
\begin{equation*} \klabel{kra 1}
\mathbf{u}={\mathbf{u}_{\mathcal{A}}}+{\mathbf{u}_{\mathbf{V}}}, \quad \text{где}\quad 
{\mathbf{u}_{\mathcal{A}}}={\lambda}^{-1}\mathbf{f}_\mathcal{A}, \quad 
\mathbf{u}_\mathbf{V}=(S+\lambda \, I)^{-1} \mathbf{f}_\mathbf{V}.
\eqno{(4.14)}\end{equation*}

Это решение $\mathbf{u}\in  \mathbf{C}(2k, m+1) $ при	$\mathbf{f}\in  \mathbf{C}(2k, m)$.
	В частности,
	
 $\mathbf{u}={\lambda}^{-1}\mathbf{f}_{\mathcal{A}}$, 	
	 если
 $\mathbf{f}_{\mathcal{A}}\in \mathcal{A}$ или  $\mathbf{f}_{\mathcal{A}}\in\mathcal{A}_{\gamma}$, а
$\mathbf{f}_{\mathbf{V}}=0$;

$\mathbf{u}=(S+\lambda  \, I)^{-1}\mathbf{f}_\mathbf{V}\in \mathbf{W}^1$
если в $\mathbf{f}\in \mathcal{B}\bot \mathcal{A}$, а $\mathbf{f}_{\mathcal{A}}=0$;

	$\mathbf{u}\in  \mathbf{C}(2k, m+1)= \mathbf{A}^{2k} \oplus \mathbf{W}^{m+1}$, 
	если 	$\mathbf{f}\in  \mathbf{C}(2k, m)$,  где   $k, m > 0 $  -целые.

Если	же $\mathbf{f}\in \mathbf{C}^{\infty}_0(\Omega)$, то поле
 $\mathbf{u}\in  \mathbf{C}^{\infty}(\overline{\Omega})$ есть классическое решение задачи.
\end{theorem}

 В случае шара,  согласно   п.1.9,  из теоремы 8 вытекает
 
 	\begin{corol} 	Если область $\Omega=B$,  %есть шар, 
 	$\psi_n(\lambda\,R)\neq 0$   $\forall\, n\in {\mathbb {N}}$,   числа $k, m$ целые, а поле  $\mathbf{f}\in  \mathbf{A}^{2k}(B)\oplus  \mathbf{W}^m(B)$,   то решение  задачи 1.1  существует, единственно и принадлежит классу    $ \mathbf{A}^{2k}(B)\oplus  \mathbf{W}^{m+1}(B)$.  	\end{corol}
 
Более того, из Теорем 5 и  8 вытекает свойство отображения $\mathrm{rot}+\lambda I$: 
 \begin{lemma}\klabel{Ro}
 	При $\lambda \overline{\in} Sp (\mathrm{rot})$ и  целом $m>0$ оператор  $\mathrm{rot}+\lambda I$  (и обратрый)
 	отображает класс  $ \mathbf{C}(2k, m+1)$
 	на  $ \mathbf{C}(2k, m)$	взаимно	однозначно и непрерывно.  \end{lemma}
 
Следующие теоремы доказываются аналогично. 

\begin{theorem}\klabel{2RG}
 При $\lambda\neq Sp(-\nabla \mathbf {div})$ единственное решение задачи 2.1  имеет вид:
\begin{equation*} \klabel{kra 1}
\mathbf{w}={\mathbf{w}_{\mathcal{A}}}+{\mathbf{w}_{\mathbf{V}}}, \quad \text{где}\quad 
{\mathbf{w}_{\mathcal{A}}}=(\mathcal{N}_d+\lambda  I)^{-1}\mathbf{f}_\mathcal{A}, \quad 
\mathbf{w}_\mathbf{V}=\lambda^{-1} \mathbf{f}_\mathbf{V}.
\eqno{(4.15)}\end{equation*}

Это решение $\mathbf{w}\in  \mathbf{C}(2(k+1), m) $ при	$\mathbf{f}\in  \mathbf{C}(2k, m)$.

	В частности, 
$\mathbf{w}=(\mathcal{N}_d+\lambda  I)^{-1}\mathbf{f}_\mathcal{A}\in \mathbf{A}^{2}(\Omega)$  при  $\mathbf{f}_\mathcal{A}\in \mathcal{A}_{\gamma}$,  $ \mathbf{f}_\mathbf{V}=0$;

$\mathbf{w}=\lambda^{-1} \mathbf{f}_\mathbf{V}$	 при  $\mathbf{f}_\mathcal{A}=0$,
 $\mathbf{f}_{\mathbf{V}}\in \mathcal{B}$;

$\mathbf{w}\in   \mathbf{C}(2(k+1), m)= \mathbf{A}^{2(k+1)}(\Omega) \oplus \mathbf{W}^{m} (\Omega) $ при	$\mathbf{f}\in  \mathbf{C}(2k, m)$. 	
\end{theorem}

 В случае шара,  согласно   п.1.9,  из теоремы 9 вытекает

\begin{corol} Если область $\Omega=B$,
	$\psi'_n(\nu\,R)\neq 0$ $\forall\, n \geq 0$,   $k, m$  целые,  а поле  $\mathbf{f}\in  \mathbf{A}^{2k}(B)\oplus  \mathbf{W}^m(B)$, то решение задачи 2  при  $p=1$ и $\lambda=\nu^2R^2$ существует, единственно и принадлежит классу    $ \mathbf{A}^{2(k+1)}(B)\oplus  \mathbf{W}^{m}(B)$.
\end{corol}

 Более того,  из Теорем 3  и 9 вытекает свойство отображения  $ \nabla \,  \mathrm{div}+\nu^2\,I$:
\begin{lemma}\klabel{o}   При     $\nu^2\neq Sp(-\nabla \, \mathrm{div})$,   
	$k > 0 $,  оператор  $ \nabla \,  \mathrm{div}+\nu^2\,I$  отображает  $ \mathbf{C}(2(k+1), m)$ на класс $ \mathbf{C}(2k,m)$ взаимно однозначно и непрерывно. \end{lemma}
Эти утверждения говорят о хорошем соответствии пространств и операторов.

\begin{theorem}\klabel{3RG}
 При $\lambda\neq Sp (-\nabla \mathrm{div})\cup Sp (rot)$ единственное решение задачи 3 имеет вид:
\begin{equation*} \klabel{kra 1}
\mathbf{u}={\mathbf{u}_{\mathcal{A}}}+{\mathbf{u}_{\mathbf{V}}}, \quad \text{где}\quad 
{\mathbf{u}_{\mathcal{A}}}=(\mathcal{N}_d+\lambda  I)^{-1}\mathbf{f}_\mathcal{A}, \quad 
\mathbf{u}_\mathbf{V}=(S+\lambda \, I)^{-1} \mathbf{f}_\mathbf{V}.
\eqno{(4.16)}\end{equation*}

Это решение $\mathbf{u}\in  \mathbf{C}(2(k+1), m+1) $ при	$\mathbf{f}\in  \mathbf{C}(2k, m)$,  где $k, m  $ целые.
\end{theorem}

{\it Замечание.}
{\small Используя теоремы 4 и 5,  читатель может продолжить 	
	 список модельных  задач и их решений.
	 Мы  предполагаем также вернуться к задаче   \cite{saUMJ11}. } 
%приложениях в (см. cdtgt,chake,motva} и  п.1.11)сообщении\quad  
 \section{Поля Бельтрами в  астрофизике и  физике плазмы}
  В  докладе НАН  \cite{chawo} \  "О бессиловых магнитных полях" \ S. Chandrasekhar and L.Woltjer отмечают, что   поля  (в областях с низкой плотностью)  являются бессиловыми в том смысле, что сила Лоренца 
 $L = j \times  H$  в них обращается в нуль ($ j$ - плотность тока, $ H$ - напряженность поля). Например, в звездных оболочках, газовых туманностях, межзвездном пространстве.
 Тот факт, что магнитное поле в Крабовидной туманности имеет необычайно регулярную структуру, значительно укрепило эту точку зрения.

 Условие $L = 0$,  очевидно,  требует, чтобы ток $ j$  протекал всюду параллельно магнитному полю; и, поскольку $4 \pi  j = curl H$, эквивалентное выражение это
 го  условия имеет вид $ curl H = \alpha H$,  где $\alpha$ - некоторая функция положения,
 произвольная, за исключением требования $\nabla \alpha \cdot  H = 0$,
 которое  вытекает из расходимости  предыдущего уравнения.
 
 Аргументы, которые обычно выдвигаются в поддержку мнения о том, что космические магнитные поля могут в значительной степени быть свободными от сил, следующие: ионизированный звездный материал имеет такую высокую электропроводность, что могут течь большие токи... 
 
 Вот, что   пишут  H.Qin, W.Liu, H.Li, and J.Squire в отчёте  PPPL-4823  \footnote{ Prinston Plasma Physics Laboratory  "Woltjer-Taylor state without Taylor's conjecture-plasma relaxation  at all wavelengths"\, Okt. 2012. } {\it
 	о релаксации плазмы:}
 "в астрофизике  и в  физике  рлазмы было открыто, что плазма релаксирует к  состоянию Woltjier- Taylor'а  \cite{wol},  \cite{tay},  которое определяется уравнениями  Бельтрами $\nabla\, \times\,\textbf{B}= \alpha\, \textbf{B}$ с постоянной $\alpha$"\  и добавляют:
 
"чтобы объяснить, как достигается такое расслабление,
Тейлор разработал свою знаменитую теорию релаксации.
основанную на предположении, 
что в релаксации преобладают коротковолновые флуктуации.

Однако не существует убедительных экспериментальных и численных доказательств, подтверждающих гипотезу Тейлора.

Разработана новая теория, которая предсказывает, что система будет развиваться в направлении W-T состояния для произвольной  флуктуации
спектра." (Конец цитаты.)

 Как мы писали ранее  в \cite{saVS20},   собственные поля ротора   называются полями Бельтрами  в гидродинамике  \cite{lad}; в астрофизике и физике плазмы-- бессиловыми полями (force-free magnetic 	fields - L.Woltjel \cite{wol},  free-decay fields -Taylor \cite{tay}).
Отметим ещё  препринт P.D.Mininni and D.C. Montgomery \cite{mimo} о численном методе решения уравнений магнитогидродинамики в сферической геометрии, а также работы 	В.И.Арнольда 
по проблеме {\it динамо}  в юбилейном сборнике \cite{ar} 1997 и
 книгу В.В.Козлова "Общая теория вихрей"\  \cite{koz} 1983, в которой имеется исторический обзор теории.

{\bf Благодарности: } Академику РАН, профессору В.П.Маслову,  профессору, доктору ф.-м. наук М.Д.Рамазанову, 
профессору,	доктору ф.-м. наук С.Ю.Доброхотову, 	доктору ф.-м. наук 
Б.И.Сулейманову и  кандидату ф.-м. наук Р.Н.Гарифуллину за  поддержку.

\end{fulltext}
%\newpage 

***

\newpage

Реферат:

Операторы вихрь и градиент дивергенции в  пространствах  
$ \mathbf{W}^{m}$ и  $\mathbf{A}^{2k}$    вихревых и потенциальных полей   и  в классах  $\mathbf{C}(2k, m)$

 Р.С.Сакс
 
{\small   Изучаются свойства  операторов вихрь    и градиент дивергенции %($ \text{rot}$ и	$\nabla \text{div}$)
	в   пространстве $\mathbf{L}_{2}(G)$ 	
	в ограниченной области   $ G\subset  \textrm{R}^3$ с гладкой границей $\Gamma$,     в	подпространствах  	$ \mathbf{W}^{m}$ и  $\mathbf{A}^{2k}$ и 	
	 в  классах    $ \mathbf{C}(2k, m)\equiv \mathbf{A}^{2k}(G) \oplus \mathbf{W}^m(G)$.  
	
		$\mathbf{L}_{2}(G) $  разлагается на ортогональные подпространства  $\mathcal{A}$  и $\mathcal{B}$:	$\mathbf{L}_{2}(G)=\mathcal{A}\oplus \mathcal{B}$. 	\newline
	Далее  
	$ \mathcal{A}= \mathcal{A}_H\oplus \mathbf{A}^0$  и  $ \mathcal{B}=\mathcal{B}_H \oplus \mathbf{V}^0$,
	где 	$\mathcal{A}_H $ и $\mathcal{B}_H $- нуль-пространства операторов 	$\nabla \text{div}$   и  $ \text{rot}$  в $\mathcal{A}$  и в  $\mathcal{B}$; 	
	размерности   $\mathcal{A}_H $ и $\mathcal{B}_H $ конечны и определяется топологией границы,  $\mathcal{A}_H=\emptyset $ и $\mathcal{B}_H= \emptyset $ в шаре  $\Omega=B $.

	Собственные поля	
	оператора	$\nabla \text{div}$  ( и $ \text{rot}$)  с ненулевыми  собственными значениями  % $\mu_{j}\neq 0$ служат основой 
	используются при построении  ортонормированного базиса  в  $ \mathbf{A}^0$  (и  в $\mathbf{V}^0$, соотв.-но). 	 
		
	Операторы   $\nabla\mathrm{div}$ и  $\mathrm{rot}$ аннулируют друг друга	  и 	проектируют $\mathbf{L}_{2}(G) $   на   $\mathcal{A}$  и $\mathcal{B}$,  причем $\mathrm{rot}\, \mathbf {u}=0$ при $\mathbf{u}\in	\mathcal{A}$,  а $\nabla\mathrm div \mathbf {v}=0$ при $\mathbf{v}\in	\mathcal{B}$.  %  \cite{hw}.	
	
	Оператор Лапласа в  $\mathbf{L}_{2}(G)$
	выражается через них: 
	$\mathrm{\Delta} \mathbf {v}
	\equiv \nabla \mathrm{div}\,\mathbf {v}
	-(\mathrm{rot})^2\, \mathbf {v}$.	
	
	Поэтому 
	$\Delta^m \mathbf{v}\equiv  (\nabla \mathrm{div})^m\mathbf{v}$,  если    $  \mathbf{v}\in \mathcal{A}$,  и  $ \Delta^m \mathbf{u} \equiv(-1)^m(\mathrm{rot})^{2m}\, \mathbf{u}$,   если       $\mathbf{u}\in \mathcal{B}$,           
	
	Кроме того, шкалы  	$ \{\mathbf{A}^{2k}\}$  и $\{\mathbf{W}^m\}$  пространств 	$ \mathbf{A}^{2k}(G)$ и $\mathbf{W}^m(G)$ представляют   	шкалу  $\{\mathbf{H}^{n}\}$  пространств Соболева  $\mathbf{H}^{n}(G)$  в ортогональных подпространствах   $\mathbf{A}^0$ и  $ \mathbf{V}^0 $ в $\mathbf{L}_{2}(G)$,   а именно,
если поле  $\mathbf{f}\in \mathbf{H}^n(G)$ и  $n\geq 2k, n\geq m$, то его проекции   $\mathbf{f}_\mathcal{A}$ и  $\mathbf{f}_\mathbf{V}$ на подпространства 
$\mathbf{A}^0$ и  $ \mathbf{V}^0 $   принадлежат  $\mathbf{A}^{2k}(G)$  и $\mathbf{W}^m(G) $, соотв.-но.  

Обратно, 
$\mathbf{A}^{2k}(G) \subset \mathbf{H}^{2k}(G)$, \quad $\mathbf{W}^m(G) \subset 
\mathbf{H}^{m}(G)$ при   $k\geq 1, m\geq1$,  значит, 
шкалы $ \{\mathbf{A}^{2k}\}$  и $\{\mathbf{W}^m\}$ 
играют роль    	шкалы   $\{\mathbf{H}^{n}\}$  в    $\mathbf{A}^0$ и  $ \mathbf{V}^0 $. 
}
В пространствах 	$ \mathbf{A}^{2k}(G)$ и $\mathbf{W}^m(G)$ 
действуют операторы	 $\mathcal{N}_d$ и $S$, -  самосопряженные расширения 
$\nabla \text{div}$ и  $ \text{rot}$    в  $ \mathbf{A}^0$ и  $\mathbf{V}^{0}$.  Они осущесвляют гомеоморизмы соседних пространств, а их сгепени   $\mathcal{N}^{2k}_d$ и $S^{2m}$- гомеоморизмы  пространств $ \mathbf{A}^{2k}(G)$ и $ \mathbf{A}^{-2k}(G)$ и  $\mathbf{W}^m(G)$  и  $\mathbf{W}^{-m}(G)$, соответственно.

Рассмотрены краевые задачи  для   $\nabla \text{div}+\lambda I$ 	 	   и	 	$ \text{rot}+\lambda I$   и других модельных операторов в  классах $  \mathbf{C}(2k, m)$. В частности, доказаны

Л е м м а 4.1.
{\it	При $\lambda \overline{\in} Sp (\mathrm{rot})$ и  целом $m>0$ операторы  $\mathrm{rot}+\lambda I$  (и обратрый)
отображают класс  $ \mathbf{C}(2k, m+1)$
на  $ \mathbf{C}(2k, m)$	взаимно	однозначно и непрерывно. }

Л е м м а 4.2.
{\it	При $ \nu ^2 \neq Sp (- \nabla \,
	\mathrm {div}) $, и  целом $k>0$ операторы   $ \nabla \, \mathrm{ div} + \nu^2 \, I $ (и обратрый)
	отображают класс  $ \mathbf{C}(2 (k + 1), m)$
	на  $ \mathbf{C}(2k, m)$	взаимно	однозначно и непрерывно. }

	Приведены аргументы астрофизиков и специалистов по физике плазмы  в пользу существования безсиловых  полей 
(полей Бельтрами) в Природе.

В статье продолжено исследование  свойств пространств	С.Л.Соболева  \cite{sob,sob54}.

*

\begin{keywords}  пространство Лебега и пространства  Соболева, операторы градиент, дивергенция,  ротор, потенциальные и вихревые поля, поля Бельтрами,	эллиптические  краевые и спектральные задачи.\end{keywords}

*

\newpage 
	  
Abstract: Saks Romen Semenovich

Сurl and gradient of divergence operators
 in Spaces $ \mathbf{W}^{m}$ and $\mathbf{A}^{2k}$ vortex and potential fields and in the classes $\mathbf{C}(2k, m)$.

The properties of the curl and the  gradient of divergence  operators 
( $ \text{rot}$	 and  $\nabla\text{div}$  ) are studied
in the space $ \mathbf {L} _{2} (G) $ 
in a bounded domain $ G \subset \textrm {R}^3 $ with a smooth boundar
 $ \Gamma$ and in the classes  $ \mathbf{C}(2k, m)(G)\equiv \mathbf{A}^{2k}(G) \oplus \mathbf{W}^m(G)$. 

The space $ \mathbf {L}_{2} (G) $ is decomposed into orthogonal subspaces   $ \mathcal{A} $ and $ \mathcal {B} $:
	$\mathbf{L}_{2}(G)=\mathcal{A}\oplus \mathcal{B}$.	
	In turn,  $ \mathcal{A}= \mathcal{A}_H\oplus \mathbf{A}^0$ and  $\mathcal{B}=\mathcal{B}_H \oplus \mathbf{V}^0$,
where  $\mathcal{A}_H $ и $\mathcal{B}_H $ are null spaces of operators
  	$\nabla \text{div}$ and  $ \text{rot}$  in $\mathcal{A}$  and $\mathcal{B}$;   the dimensions of   $\mathcal{A}_H $ and $\mathcal{B}_H $   are finite and  determined by the topology of the boundary;   $\mathcal{A}_H=\emptyset $ and $\mathcal{B}_H= \emptyset $  if  the domain $\Omega $ is a ball.

  	 The orthonormal basis are constructed in the class  $ \mathbf{A}^0$ (resp.,  In $\mathbf{V}^0$ ) 
  	  by   	   eigenfields
  	$\mathbf{q}_{j}(\mathbf{x})$  of	$\nabla \text{div}$ operatoг (resp.,	$\mathbf{q}^{\pm }_{j}(\mathbf{x})$ of	$ \text{rot}$  operatoг) with nonzero eigenvalues    $\mu_{j}$   (resp., $\pm \lambda_{j}$ ).

  The operators $\nabla\mathrm{div}$ and $\mathrm{rot}$ cancel each other out and project $\mathbf{L}_{2}(G) $  onto $ \mathcal {A} $ and $ \mathcal { B} $, and $ \mathrm {rot} \, \mathbf {u} = 0 $ for $ \mathbf {u} \in \mathcal {A} $, and $ \nabla \mathrm div \mathbf {v} = 0 $ for $ \mathbf {v} \in \mathcal {B} $ \cite{hw}.
  	
Laplace matrix operator expressed through them:
 $\mathrm{\Delta} \mathbf {v}
 \equiv \nabla \mathrm{div}\,\mathbf {v}
 -(\mathrm{rot})^2\, \mathbf {v}$.	
 
 Therefore
 $\Delta^m \mathbf{v}\equiv  (\nabla \mathrm{div})^m\mathbf{v}$,  if    $  \mathbf{v}\in \mathcal{A}$,  and  $ \Delta^m \mathbf{u} \equiv(-1)^m(\mathrm{rot})^{2m}\, \mathbf{u}$,   if     $\mathbf{u}\in \mathcal{B}$.
  
  Moreover, the scales $ \{\mathbf{A}^{2k}\}$ and $\{\mathbf{W}^m\}$ of the spaces $ \mathbf{A}^{2k}(G)$ and $ \mathbf{W}^m(G)$ represent the scale of $\{\mathbf{H}^{n}\}$ Sobolev spaces $\mathbf{H}^{n}(G)$ in the orthogonal subspaces of $\mathbf {A}^0$ and $ \mathbf{V}^0 $ in $\mathbf{L}_{2}(G)$, namely,
  
if the field  $\mathbf{f}$   belong to  $ \mathbf{H}^n(G)$ and $n\geq 2k, n\geq m$, then  $\mathbf{f}_\mathcal{A}$ and $\mathbf{f}_\mathbf{V}$  (  its projections into subspaces
$\mathbf{A}^0$ and $ \mathbf{V}^0 $)  belong to $\mathbf{A}^{2k}(G)$ and $\mathbf{W}^m(G) $, respectively. 

Back,
$\mathbf{A}^{2k}(G) \subset \mathbf{H}^{2k}(G)$, \quad $\mathbf{W}^m(G) \subset
\mathbf{H}^{m}(G)$ for $k\geq 1, m\geq1$, then
scales $ \{\mathbf{A}^{2k}\}$ and $\{\mathbf{W}^m\}$
play the role of the scale $\{\mathbf{H}^{n}\}$ in $\mathbf{A}^0$ and $ \mathbf{V}^0 $.

The operators $\mathcal{N}_d$ and $S$ act
in the spaces $ \mathbf{A}^{2k}(G)$ and $\mathbf{W}^m(G)$, they are self-adjoint extensions
$\nabla \text{div}$ and $ \text{rot}$ into $ \mathbf{A}^0$ and $\mathbf{V}^{0}$.

They realize homeomorisms of neighboring spaces, and their degrees $\mathcal{N}^{2k}_d$ and $S^{2m}$ are homeomorisms of the spaces $ \mathbf{A}^{2k}(G)$ and $ \mathbf {A}^{-2k}(G)$ and $\mathbf{W}^m(G)$ and $\mathbf{W}^{-m}(G)$, respectively.

Boundary value problems for $\nabla \text{div}+\lambda I$ and $ \text{rot}+\lambda I$ and other model operators in the classes $ \mathbf{C}(2k, m)$ are considered.
In particular, it has been proven

{\it  	For $ \lambda \overline {\in} Sp (\mathrm {rot}) $ and
	$  m> 0$, the operators  $ \mathrm {rot} + \lambda I $ (and its inverse) maps the class $ C (2k, m + 1) $
	on $ C (2k, m) $  one-to-one and continuous} (Lemma 4.1). 

{\it  	 For $ \nu ^2 \neq Sp (- \nabla \,
	\mathrm {div}) $, $ k> 0 $, operators $ \nabla \, \mathrm{ div} + \nu^2 \, I $ maps $ C (2 (k + 1), m) $ to the class $ C (2k, m) $ one-to-one and continuous} (Lemma 4.2)

The arguments of astrophysicists and specialists in plasma physics in favor of the existence of force-free fields (that is, Beltrami fields) in nature are presented.

The paper continues the study of the properties of Sobolev spaces \cite{sob,sob54}.

*

{\bf Keywords:} Lebesgue space and Sobolev spaces, operators gradient, divergence, rotor (curle), potential and vortex fields,   Beltrami fields, elliptic boundary value and spectral problems.

*

Сакс Ромэн Семенович

старший научный сотрудник

 Институт Математики с ВЦ УФИЦ РАН 

450077, г. Уфа, ул. Чернышевского, д.112 

телефон: (347)272-59-36 \quad  (347) 273-34-12

факс: (347) 272-59-36

 телефон дом.: (347) 273-84-69
 
  моб.+7 917 379 75 38

 e-mail: romen-saks@yandex.ru

*

%\end{fulltext}

\end{document}